\documentclass[12pt]{amsart}
\usepackage{latexsym,enumerate}
\usepackage{amsmath,amsthm,amsopn,amstext,amscd,amsfonts,amssymb,mathrsfs}
\usepackage{color}
\usepackage{color}
\numberwithin{equation}{section}

\newtheorem{theorem}{Theorem}
\newtheorem{lemma}{Lemma}
\newtheorem{corollary}{Corollary}
\newtheorem{proposition}{Proposition}
\newtheorem{remark}{Remark}

\numberwithin{theorem}{section}
\numberwithin{corollary}{section}
\numberwithin{lemma}{section}
\numberwithin{definition}{section}
\numberwithin{proposition}{section}
\numberwithin{remark}{section}

\newcommand{\R}{\mathbb R}

\newcommand{\medint}{-\kern  -,375cm\int}
\newcommand{\dint}{\displaystyle\int}

\begin{document}

\title[]{Some isoperimetric inequalities on $\mathbb{R} ^N$
\\
with respect to weights $|x|^\alpha $}

\author[A. Alvino]{A. Alvino$^1$}
\author[F. Brock]{F. Brock$^2$}
\author[F. Chiacchio]{F. Chiacchio$^1$}
\author[A. Mercaldo]{A. Mercaldo$^1$}
\author[M.R. Posteraro]{M.R. Posteraro$^1$}

%\date{\today}

\setcounter{footnote}{1}
\footnotetext{Dipartimento di Matematica e Applicazioni ``R. Caccioppoli'',
Universit\`a degli Studi di Napoli Federico II,
Complesso Monte S. Angelo, via Cintia, 80126 Napoli, Italy;\\
e-mail: {\tt angelo.alvino@unina.it, fchiacch@unina.it,  mercaldo@unina.it, posterar@unina.it}}

\setcounter{footnote}{2}
\footnotetext{University of Rostock, Department of 
Mathematics, 
Ulmenstr. 69, 18057 Rostock, Germany, email: friedemann.brock@uni-rostock.de}

\begin{abstract} 
We solve a class of isoperimetric problems on $\mathbb{R}^N $ with respect to weights that are powers of the distance to the origin. For instance  we show that, if $k\in [0,1]$, then  among all smooth sets $\Omega$ in $\mathbb{R} ^N$ with fixed Lebesgue measure, $\int_{\partial \Omega } |x|^k \, \mathscr{H} _{N-1} (dx)$ achieves its minimum for a ball centered at the origin. Our results also imply a weighted Polya-Sz\"ego principle. In turn,  we establish radiality of optimizers  in some
Caffarelli-Kohn-Nirenberg inequalities, and  we obtain  
 sharp bounds for eigenvalues of some nonlinear problems.

\medskip

\noindent
{\sl Key words: isoperimetric inequality, weighted rearrangement, norm inequality, elliptic boundary value problem, eigenvalue problem}  
\rm 
\\[0.1cm]
{\sl 2000 Mathematics Subject Classification:} 51M16, 46E35, 46E30, 35P15 
\rm 
\end{abstract}
\maketitle

\section{Introduction}

There has been a growing interest in isoperimetric inequalities with weights during the last decades and a wide literature is available, see for instance, \cite{BBMP}, \cite{BBMP2}, \cite{BBCLT}, \cite{BMP},   \cite{XR},   \cite{XRS}, \cite{CMV},   \cite{CJQW}, \cite{Cham},   \cite{C},\cite{DDNT}, \cite{DHHT}, \cite{KZ}, \cite{Mo}, \cite{Mo2}  and  the  references therein. However, until now most research dealt with inequalities where both the volume functional and perimeter functional carry the same weight.
In this article we make analyse a scale of isoperimetric inequalities on 
$\mathbb{R}^N$ with two 
{\sl different } weights in perimeter and 
volume which are powers of the distance to the origin. 
\\
More precisely, given $k,l\in \mathbb{R}$, we study the following isoperimetric problem:
\medskip

{\sl  Minimize $\displaystyle\int_{\partial \Omega } |x|^k \, \mathscr{H}_{N-1} (dx) $ among all smooth sets $\Omega \subset \mathbb{R} ^N$ satisfying $\dint_{\Omega } |x|^l \, dx =1$.}
\medskip

In particular, we are interested in conditions on the numbers $k$ and $l$ such that 
the above minimum is attained for a ball centered at the origin.
\\  
Our motivation comes from some norm inequalities with weights
which are now well-known as the {\sl Caffarelli-Kohn-Nirenberg inequalities} (see, e.g. \cite{CKN}, \cite{CW}, \cite{CP}, \cite{DolEstLoss} and the references cited therein).
These inequalities compare a weighted $L^p $- norm of the gradient of a 
function on $\mathbb{R}^N $ with a weighted $L^q $--norm 
of the function, and they have many applications to the analysis of weighted elliptic and parabolic problems.
\\ 
%Let us summarize the main symmetry results for smooth sets as they are available in the literature now. 
Let us state the main results concerning the isoperimetric inequalities. We emphasize that the result has been already proved in the cases (i) and (ii) in \cite{Howe}, respectively \cite{ChiHo}.
\begin{theorem}
\label{maintheorem}
Let $N\in \mathbb{N} $, 
$k,l\in \mathbb{R} $ and $l+N >0$. Further, assume that one of the following conditions holds:
\\
{\bf (i)} $N\geq 1 $ and $l+1\leq k $;
\\
{\bf (ii)} $N\geq 2$, $k\leq l+1$ and $ l\frac{N-1}{N} \leq k\leq 0$; 
\\  
{\bf (iii)} $N\geq 3$, $ 0\leq k\leq l+1$ and 
$$
\frac{1}{l+N } \geq \frac{1}{k+N-1} - \frac{ (N-1)^2 }{ N(k+N-1)^3 }.
$$
{\bf (iv)} $N= 2 $, $k\leq l+1$, and either 
\begin{eqnarray*}
 & & l\leq 0\leq k\leq \frac{1}{3} \ \mbox{ or }
 \\
 & & \frac{1}{3} \leq k \ \mbox{ and } \
 \frac{1}{l+2} \geq \frac{1}{k+1} - \frac{16}{27 (k+1)^3 } ;
\end{eqnarray*} 
\\
Then  
\begin{equation}
\label{mainineq}
\dint_{\partial \Omega } |x|^k \, \mathscr{H}_{N-1} (dx)
\geq 
C_{k,l,N} ^{rad}  
\left( 
\dint_{\Omega } |x|^l \, dx 
\right) 
^{(k+N-1)/(l+N) } , 
\end{equation}
for all smooth sets $\Omega $ in $\mathbb{R}^N $,
where 
\begin{equation}
\label{defCkl}
C_{k,l,N} ^{rad} := (N\omega _N ) ^{(l-k+1)/(l+N)} \cdot  (l+N ) ^{(k+N-1)/(l+N)} .
\end{equation}
Equality in (\ref{mainineq}) holds  for every ball centered at the origin.
\end{theorem}

\noindent
Let us briefly comment on Theorem \ref{maintheorem}. 
\\
First observe that Theorem \ref{maintheorem} can be extended to Lebesgue measurable sets on $\mathbb{R}^N $ by a standard approximation procedure (see Section 2). Moreover, it is often possible to detect all
cases of equality in (\ref{mainineq}) (for details, see Section 5). 
Inequality (\ref{mainineq}) was proved under assumption {\bf (i)} in \cite{Howe}, Theorem 3.1 and its application through Example 3.5, part (4), and under assumption {\bf (ii)} in  \cite{ChiHo}, Theorem 1.3. Note also that partial results in the cases {\bf (i)} and {\bf (ii)} were obtained in \cite{DHHT}, Proposition 4.21, part (2) and (3), and in \cite{BBMP}, Theorem 2.1. 
\\
The main result of this paper is the proof of  Theorem \ref{maintheorem} in the cases {\bf (iii)} and {\bf (iv)}. 
We emphasize that the conditions {\bf (ii)}-{\bf (iv)} contain the range $N\geq 2 $, $l=0\leq k\leq 1$, while the case $l=0, k\ge 1$ was already known for some time, see \cite{BBMP}. In this way we generalize in particular on a recent result in \cite{C} where only the two-dimensional case was considered.  
\\
We wish to point out that the situation can be quite different from Theorem \ref{maintheorem} for other ranges of the parameters $k$ and $l$. For instance, if $k=l\ge 0$, then the minimizing sets have been identified as balls whose boundaries touch the origin, see \cite{BBCLT} and  \cite{DDNT}. When $N=2$, then the same result even holds for the range $0\le k\le l\leq  2k$, see \cite{DHHT}, Proposition 4.21, part (4).
\\
Now we outline the content of the paper. We introduce some notation and provide some analytic tools that will be of later use in Section 2. In Section 3 we introduce two functionals ${\mathcal R}_{k,l,N} $ and ${\mathcal Q}_{k,l,N} $ and we provide some basic information related to the isoperimetric problem.  
In Section 4 we give a necessary condition for the existence of a minimizer (Lemma \ref{R3}) and a necessary condition for the radiality of minimizers in the isoperimetric problem, which also shows a symmetry beaking (see Theorem \ref{R4}). 
Section 5 deals with the proof of Theorem \ref{maintheorem}. In addition, we treat the equality case in (\ref{mainineq}), 
see the Theorems \ref{R5}, \ref{th1bis}, \ref{th1ter}, Corollary \ref{N=2_l<0} and Theorem \ref{theoN=2k>0}. Further, we give a complete solution of the isoperimetric problem in the case $N=1$ in Section 6, Theorem \ref{1d_Iso}. Our proofs use well-known rearrangement tools, the classical isoperimetric inequality and Hardy's inequality.
The interpolation argument that occurs in the proof of Lemma \ref{4.3} 
seems to be new in this context.
Using the previous results and inversion in the unit sphere, we show 
an isoperimetric inequality where the extremal sets are exteriors of balls centered at the origin in Theorem \ref{maintheorem2}. Finally, we give some applications of Theorem \ref{maintheorem} in Section 8. Using the notion of weighted rearrangement we provide a Polya-Sz\"ego-type inequality in 
Theorem \ref{ps}. Then we use this to obtain best constants in some Caffarelli-Kohn-Nirenberg inequalities (see Lemma  \ref{CKN}, Proposition \ref{improve1} and Theorems \ref{2ndimprove} and \ref{3rdimprove}). Further, in Theorem \ref{lorentz} we obtain 
the best constant in a weighted Sobolev-type inequality for Lorentz spaces, originally proved in \cite{Alvino} (see also \cite{E}, \cite{CRT}), and 
a sharp bound for the first eigenvalue of a weighted elliptic eigenvalue problem associated to the $p$-Laplace operator (Theorem \ref{eigenvalue}).

\section{Notation and preliminary results}

Throughout this article  $N$ will denote a natural number while   
$k$ and $l $ are real numbers. With the exception of Section 5 we will assume
\begin{equation}
\label{ass1}
k+N-1 >0 \ \mbox{ and } \ l+N>0 .
\end{equation}

\noindent 
Let us introduce some notation. 
\begin{eqnarray*}
B_{R}(x_0 )
 & := &
\left\{ x\in \mathbb{R}^{N}:\left\vert x-x_0 \right\vert <R\right\} , \quad (x_0 \in \mathbb{R}^N ),
\\
B_{R} 
& := & 
B_{R}(0), \quad (R>0). 
\end{eqnarray*}
$\mathscr{L} ^N $ denotes $N$-dimensional Lebesgue measure and 
$\omega_N = \mathscr{L}^N (B_1 )$.
\\
We will use frequently $N$-dimensional spherical coordinates $(r, \theta)$ in $\mathbb{R} ^N$: 
$$
\mathbb{R}^N \ni x = r\theta , \quad \mbox{where $r=|x|$, and $\theta = x|x|^{-1} \in \mathscr{S}^{N-1} $.}
$$ 
If $M$ is any set in $\mathbb{R}^N  $, then let $\chi _M $ denote its characteristic function.

Next, let $k$ and $l$ be real
numbers satisfying (\ref{ass1}). We define a measure $\mu _{l}$ by 
\begin{equation}
d\mu _{l}(x)=|x|^{l}\,dx.  
\label{dmu}
\end{equation}
If $M\subset $ ${\mathbb R}^{N}$ is a  measurable set with finite 
$\mu _l $-measure, then let 
$M^{\star }$ denote the ball $B_{R}$ such that
\begin{equation}
\mu_l \left( B_{R} \right) =\mu _l \left( M\right) = \int_M d\mu _l (x) .  
\label{mu_(M)}
\end{equation}
If $u: \mathbb{R}^N \rightarrow \mathbb{R} $ is a measurable function such that
$$
\mu_l \left( \left\{ |u(x)|>t\right\} \right) <\infty \qquad \forall t>0,
$$
then let $u^{ \star }$ denote the weighted Schwarz symmetrization of $u$, or
in short, the \\
$\mu_l -$symmetrization of $u$, which is given by
\begin{equation}
u^{ \star }(x)=\sup \left\{ t\geq 0:\mu_l \left( \left\{ |u(x)| >t\right\} \right)
>\mu _l \left( B_{\left\vert x\right\vert }\right) \right\} .  
\label{u_star}
\end{equation}

Note that $u^{\star }$ is radial and radially non-increasing, 
and if $M$ is a measurable set with finite $\mu _l $-measure, then 
$$
\left( \chi _M \right) ^{\star} = \chi _{M^{ \star }} .
$$

The {\sl $\mu _{k}$--perimeter\/} of a measurable set $M $ is given by 
\begin{equation}
P_{\mu _{k}}(M ):=\sup \left\{ \int_{M }\mbox{div}\,\left( |x|^{k}
\mathbf{v}\right) \,dx:\,\mathbf{v}\in C_{0}^{1}(\mathbb{R}^N ,\mathbb{R}^{N}),\,|
\mathbf{v}|\leq 1\mbox{ in }\, M \right\} .
\end{equation}

\noindent It is well-known that, if $\Omega $ is an open set, then the above 
\textsl{distributional definition} of
weighted perimeter is equivalent to the following

\begin{equation}
P_{\mu_{k} }(\Omega )= \left\{ 
\begin{array}{ccc}
\displaystyle\int_{\partial \Omega }|x|^{k} \, \mathscr{H} _{N-1}(dx)
 & \mbox{ if } & 
\partial \Omega  \mbox{ is } (N-1)-\mbox{rectifiable } \\ 
&  &  \\ 
+ \infty \qquad
 & \mbox{ otherwise } & 
\end{array}
\right.
.
\end{equation}
($\mathscr{H} _{N-1} $ denotes $(N-1)$-dimensional Hausdorff-measure.)  

We will call a set $\Omega \subset \mathbb{R}^N $ {\sl smooth}, 
if it is open and  bounded with smooth boundary, that is, for every $x_0 \in \partial \Omega $, 
there is a number $r >0$ such that $B_r (x_0 ) \cap \Omega $
has exactly one connected component and $B_r (x_0 ) \cap \partial \Omega $ 
is the graph of a $C^1 $--function on an open set in $\mathbb{R} ^{N-1} $.  

If $p\in \left[ 1,+\infty \right) $, we will denote by $L^{p}(\Omega ,d\mu
_{l})$ the space of all Lebesgue measurable real valued functios $u$ such
that
\begin{equation}
\left\Vert u \right\Vert 
_{L^{p}(\Omega ,d\mu _{l})}
:=\left( 
\int_{\Omega
}
\left\vert 
u\right\vert 
^{p}d\mu _{l} (x) 
\right) ^{1/p}
<+\infty .
\label{Norm_Lp}
\end{equation}
We will often use the following well-known {\sl Hardy-Littlewood inequality}, 
\begin{equation}
\label{hardylitt1}
\int_{\mathbb{R}^N } uv \, d\mu _l (x) \leq \int_{\mathbb{R}^N } u^{ \star} v^{\star} \, d\mu _l (x) ,
\end{equation}
which holds for all functions $u,v\in L^2 (\mathbb{R} ^N , d\mu _l )$.
\\  
Furthermore, $W^{1,p}(\Omega ,d\mu _{l})$ is the weighted Sobolev space
consisting of all functions which together with their weak derivatives $u_{x_{i}}$, ($i=1,...,N$), 
belong to $L^{p}(\Omega ,d\mu _{l})$
The
space will be equipped with the norm
\begin{equation}
\left\Vert u\right\Vert _{W^{1,p}(\Omega ,d\mu _{l})}:=\left\Vert
u\right\Vert _{L^{p}(\Omega ,d\mu _{l})}+\left\Vert \nabla u\right\Vert
_{L^{p}(\Omega ,d\mu _{l})}.  
\label{Norm_Wp}
\end{equation}
Finally $W_{0}^{1,p}(\Omega ,d\mu _{l})$ will stand for the closure of 
$C_{0}^{\infty }(\Omega)$ under norm (\ref{Norm_Wp}).

Now we want to recall the so-called starshaped rearrangement 
(see \cite{Kaw}) which we will use in Section 5.
For later convenience, we will write $y$ for points in $\mathbb{R}^N $  and  $(z, \theta )$ for corresponding $N$-dimensional spherical coordinates ($z= |y|$, $\theta = y|y|^{-1} $).
\\
We call a measurable set $M\subset \mathbb{R} ^N $ {\sl starshaped\/} if the set 
$$
M\cap \{ z\theta : \, z\geq 0 \}
$$ 
is either empty or a segment 
$\{ z\theta : \, 0\leq z< m(\theta ) \} $ 
for some number $m(\theta ) >0 $, for almost every $\theta \in {\mathscr S} ^{N-1} $. 
\\
If $M$ is a bounded measurable set in $\mathbb{R}^N $, and 
$\theta \in {\mathscr S}^{N-1} ,$ then let
$$
M(\theta ) := M\cap \{ z\theta :\, z\geq 0\}.  
$$
There is a unique number $m(\theta )\in [0,+\infty )$ such that
$$
\int_0 ^{m(\theta )} z^{N-1}\, dz = \int_{M(\theta )} z^{N-1} \, dz.
$$
We define 
$$
\widetilde{M}(\theta ) := \{ z\theta : \, 0\leq z\leq m(\theta )  \} ,
\quad (\theta \in {\mathscr S} ^{N-1} ),
$$
and 
$$
\widetilde{M} := \{ z\theta : \, z\in \widetilde{M}(\theta ) , \, 
\theta \in {\mathscr S} ^{N-1} \} .
$$
We call the set $\widetilde{M}$ the {\sl starshaped rearrangement of $M$\/}. 
\\
Note that $\widetilde{M} $ is Lebesgue measurable and starshaped, and 
we have
\begin{equation}
 \label{starsh1}
 {\mathscr L} ^N (M) = {\mathscr L} ^N (\widetilde{M}).
\end{equation}
If $v:\mathbb{R} ^N \to \mathbb{R} $ is measurable with compact support, 
and $t\geq 0$, then let
$ E_t $ be the super-level set $\{ y: \, |v(y)| \geq t\} $. We
define 
$$
\widetilde{v} (y) := \sup \{ t\geq 0 :\, y \in \widetilde{E_t }   \} .
$$
We call $\widetilde{v} $ the {\sl starshaped rearrangement of $v$ \/}.  
It is easy to verify that 
$\widetilde{v}$ is equimeasurable with $v$, 
that is, the following properties hold:
\begin{eqnarray}
\label{starsh2}
 & & \widetilde{E_t} = \{ y:\, \widetilde{v} (y)\geq t\} , 
 \\
 \label{starsh3}
 & & {\mathscr L} ^N (E_t ) = 
 {\mathscr L} ^N (\widetilde{E_t} ) \quad \forall t\geq 0.
\end{eqnarray}
This also implies  Cavalieri's principle: 
If $F\in C ([0, +\infty ))$ with $F(0)=0$ and if 
$F(v) \in L^1 ( \mathbb{R}^N ) $, then
\begin{equation}
\label{caval1}
\int_{\mathbb{R} ^N } F(v)\, dy =  
\int_{\mathbb{R} ^N } F(\widetilde{v} )\, dy
\end{equation}
and if $F$ is non-decreasing, then 
\begin{equation}
\label{monrearr}
\widetilde{F(v)} = F(\widetilde{v}).
\end{equation}
Note that the mapping
$$
z\longmapsto \widetilde{v} (z\theta ) , \quad (z\geq 0),
$$
is non-increasing for all $\theta \in {\mathscr S} ^{N-1} $.
\\ 
If $v,w\in L^2 (\mathbb{R}^N )$ are functions with compact support,
then there holds 
Hardy-Littlewood's inequality:
\begin{equation}
\label{harlit}
\int_{\mathbb{R}^N } vw \, dy \leq  \int_{\mathbb{R}^N } 
\widetilde{v} \widetilde{w} \, dy.
\end{equation}    
If $f:(0,+\infty) \to \mathbb{R} $ is a measurable function with compact support, then its (equimeasurable) {\sl non-increasing rearrangement }, $\widehat{f} : (0,+\infty )\to [0,+\infty )$, 
is the monotone non-increasing 
function such that 
$$
\mathscr{L} ^1 \{ t  \in [0,+\infty )  :\,  |f(t )| > c\}  = 
\mathscr{L}^1 \{ t \in [0,+\infty ) :\,  \widehat{f}(t ) > c \}  \quad \forall c\geq 0,
$$
see \cite{Kaw}, Chapter 2.
A general Polya-Szeg\"o principle for non-increasing rearrangement has been given in \cite{Lan}, Theorem 2.1. For later reference we will only need a special case: 
\begin{lemma}
\label{Landes} 
Let $\delta \geq 0$,  and let $f:(0,+\infty ) \to \mathbb{R} $ be a bounded, locally Lipschitz continuous function with bounded support, such that 
$$
\int_0 ^{+\infty } t  ^{\delta } |f' (t ) |\, dt <+\infty .
$$
Then $\widehat{f} $ is locally Lipschitz continuous and 
\begin{equation}
\label{landes1}
\int_0 ^{+\infty } t^{\delta } |\widehat{f}' (t ) |\, dt
\leq
\int_0 ^{+\infty } t^{\delta } |f' (t ) |\, d t.
\end{equation}  
\end{lemma}
 \bigskip

\section{The functionals ${\mathcal R}_{k,l,N}$ and ${\mathcal Q}_{k,l,N}$ }
Throughout this section  we  assume (\ref{ass1}), i.e.
\begin{equation*}
k+N-1 >0 \ \mbox{ and } \ l+N>0 .
\end{equation*}

If $M $ is a measurable set with $0<\mu _l (M)<+\infty $, we set
\begin{equation}
\label{rayl1}
{\mathcal R}_{k,l,N} (M) := \frac{P_{\mu _k } (M) }{ 
\left( \mu_l (M) \right) ^{(k+N-1)/(l+N) }}. 
\end{equation}
Note that
\begin{equation}
\label{Rklsmooth}
{\mathcal R} _{k,l,N} (\Omega ) = 
\frac{
\dint_{\Omega } |x|^k \, dx 
}{
\left( \dint_{\Omega } |x|^l \, dx \right) ^{(k+N-1)/(l+N)} } 
\end{equation}
for every smooth set $\Omega \subset \mathbb{R}^N $.

If $u\in C_0 ^1  (\mathbb{R} ^N )\setminus \{ 0\} $, we set
\begin{equation}
\label{rayl2}
{\mathcal Q}_{k,l,N} (u ) := \frac{\dint_{\mathbb{R} ^N } |x|^k |\nabla u| \, dx}{ 
\left( \dint_{\mathbb{R} ^N } |x|^l |u| ^{(l+N)/(k+N-1)} \, dx \right) ^{(k+N-1)/(l+N)}}. 
\end{equation}

Finally, we define 
\begin{equation}
\label{isopco}
C_{k,l,N}^{rad} := (N\omega _N ) ^{(l-k+1)/(l+N)} \cdot  (l+N ) ^{(k+N-1)/(l+N)}
.
\end{equation} 
We study the following isoperimetric problem: 
\\[0.5cm]
{\sl Find the constant $C_{k,l,N} \in [0, + \infty )$, such that}
\begin{equation}
\label{isopproblem}
C _{k,l,N} := \inf \{ {\mathcal R}_{k,l,N} (M):\, 
\mbox{{\sl $M$ is measurable with $0<\mu _l (M) <+\infty $.}} \} 
\end{equation}   
Moreover, we are interested in conditions on $k$ and $l$ such that
\begin{equation}
\label{isoradial}
{\mathcal R}_{k,l,N} (M) \geq {\mathcal R}_{k,l,N} (M^{ \star} ) 
\end{equation}
holds for all measurable sets $M$ with  $ 0<\mu _l (M)<+\infty $. 
\\[0.1cm]
Let us begin with some immediate observations.
\\
If $M$ is a measurable set with finite $\mu _l $-measure and with 
finite $\mu_k $-perimeter, then there exists a sequence of smooth sets 
$\{ M_n \} $ such that $\lim_{n\to \infty } \mu _l (M_n \Delta M) =0$ and 
$\lim_{n\to \infty } P_{\mu _k } (M_n ) = P_{\mu _k } (M)$. 
This property is well-known for Lebesgue measure (see for instance 
\cite{G}, Theorem 1.24) 
and its proof carries over to the weighted case. This implies that we also have 
\begin{equation}
\label{CklNsmooth}
C_{k,l,N} = \inf \{ {\mathcal R}_{k,l,N} (\Omega ):\, \Omega \subset \mathbb{R} ^N, \, \Omega  
\mbox{ smooth} \} .
\end{equation}  
The functionals ${\mathcal R}_{k,l,N} $ and ${\mathcal Q}_{k,l,N} $ 
have the following homogeneity properties, 
\begin{eqnarray}
\label{hom1}
 {\mathcal R}_{k,l,N} (M ) & = & {\mathcal R}_{k,l,N} (tM ) ,
\\
{\mathcal Q}_{k,l,N} (u) & = & {\mathcal Q}_{k,l,N} (u^t ),
\end{eqnarray}
where  $t>0$, $M $ is a measurable set with $0<\mu_l (M)<+\infty $, 
$u\in C_0 ^1 (\mathbb{R}^N )\setminus \{ 0\}$,  \\
$tM := \{tx:\, x\in M \} $ 
and $u^t (x):= u(tx) $, ($x\in \mathbb{R} ^N $), and there holds   
\begin{equation}
\label{isopconst2}
C_{k,l,N} ^{rad} = {\mathcal R}_{k,l,N} (B_1 ).
\end{equation}  
Hence we have that 
\begin{equation}
\label{relCC}
C_{k,l,N} \leq C_{k,l,N} ^{rad} ,
\end{equation}  
and (\ref{isoradial}) holds if and only if 
$$
C_{k,l,N} = C_{k,l,N} ^{rad} .
$$  
Finally, the classical isoperimetric inequality reads as
\begin{eqnarray}  
\label{isopclass}
 & & {\mathcal R} _{0,0,N} (M ) \geq C_{0,0,N} ^{rad} \quad 
\mbox{{\sl for all measurable sets $M$ with $0< \mu _0  (M)<+\infty $,}}
\\
\nonumber
 & & \mbox{{\sl and equality holds only if $M$ is a  ball in $\mathbb{R}^N $.}} 
\end{eqnarray} 
\begin{lemma}
\label{hardylitt}
Let $l>l' >-N $. Then 
\begin{equation} 
\label{hardylitt2}
\frac{\left( \mu _l (M)  \right) 
^{1/(l+N)} 
}{ 
\left( \mu _{l' } (M)  \right)
^{1/(l'+N)}
} 
\geq \omega_N ^{\frac{1}{l+N} -\frac{1}{l'+N} } 
\cdot 
(l+N ) ^{-\frac{1}{l+N} } 
(l' + N) ^{\frac{1}{l'+N} } 
\end{equation}
for all measurable sets $M$ with  $0<\mu_l (M)<+\infty $.
Equality holds only for balls $B_R $, ($R>0$).
\end{lemma}

{\sl Proof: } Let $M^{ \star} $ be the 
$\mu_l $-symmetrization of $M$. Then we obtain, using the Hardy-Littlewood inequality,
\begin{eqnarray*}
\mu _{l' } (M) =\int_M |x| ^{l'} \, dx & = & \int_{\mathbb{R} ^N }  |x|^{l'-l} \chi _M (x)\, d\mu _l (x) 
\\
 & \leq & 
\int_{\mathbb{R} ^N }  \left( |x|^{l'-l} \right) ^{ \star}  \left( \chi _M  \right) ^{ \star} (x)\, d\mu _l (x) 
\\
 & = & 
\int_{\mathbb{R} ^N }  |x|^{l'-l}  \chi _{M^{ \star} }  (x)\, d\mu _l (x) 
\\
 & = & 
\int_{M^{ \star} } |x|^{l' }\, dx =\mu _{l' } (M^{ \star} ).
\end{eqnarray*} 
This implies that
$$
\frac{\left( \mu _l (M)  \right) 
^{1/(l+N)} 
}{ 
\left( \mu _{l' } (M)  \right)
^{1/(l'+N)}
} 
\geq
\frac{\left( \mu _l (M^{ \star} )  \right) 
^{1/(l+N)} 
}{ 
\left( \mu _{l' } (M^{ \star} )  \right)
^{1/(l'+N)}
} 
$$
and, by  evaluating the right-hand side, (\ref{hardylitt2}) follows. 

\noindent Next assume that equality holds in (\ref{hardylitt2}). Then we must have 
$$
\int_M |x|^{l'-l} \, d\mu _l (x) = \int_{M^{ \star} } |x|^{l'-l} d\mu _l (x) ,
$$
that is, 
$$
\int_{M\setminus M^{ \star} } |x|^{l'-l} \, d\mu _l (x) = \int_{M^{ \star} \setminus M} |x|^{l'-l} d\mu _l (x) .
$$
Since $l'-l<0$, this means that 
$ \mu _l ( M\Delta M^{ \star} )=0$. The Lemma is proved. 
$\hfill \Box $
\begin{lemma}
\label{rangekl1}
Let $k,l$ satisfy (\ref{ass1}). Assume that $l>l' >-N$ and 
$C_{k,l,N}  = C_{k,l,N} ^{rad} $. 
Then we also have 
$C_{k,l',N}  = C_{k,l',N} ^{rad} $.
Moreover, if  $
{\mathcal R}_{k,l',N} (M ) = C_{k,l',N} ^{rad} $  
for some measurable set $M$ with $0< \mu _{l' } (M) <+\infty $,
then $M $ is a ball centered at the origin.
\end{lemma}
{\sl Proof:} By our assumptions and Lemma \ref{hardylitt}  we have for every measurable set $M$ with 
$0<\mu _l (M) <+\infty $,  
\begin{eqnarray*}
{\mathcal R}_{k,l',N} (M ) & = & {\mathcal R}_{k,l,N} (M ) 
\cdot 
\left[
\frac{
\left( 
\mu _l (M) 
\right) ^{1/(l+N)}
}{ 
\left( \mu_{l' } (M) 
\right) ^{1/(l'+N)}
} 
\right] ^{k+N-1}
\\
 & \geq & 
C_{k,l,N}^{rad} 
\cdot 
\left[ 
\omega_N ^{\frac{1}{l+N} -\frac{1}{l'+N} } 
\cdot 
(l+N ) ^{-\frac{1}{l+N} } 
(l' + N) ^{\frac{1}{l'+N} } 
\right] ^{k+N-1} 
= C_{k,l',N}^{rad},
\end{eqnarray*}
with equality only if 
$M $ is a ball centered at the origin. 
$\hfill \Box $

\begin{lemma}
\label{R2} 

Assume that $k \leq l+1$. Then
\begin{equation}
\label{ineqQR}
C_{k,l,N} 
=
\inf \left\{ {\mathcal Q}_{k,l,N} (u) :\, u\in C_0 ^1 (\mathbb{R} ^N 
 )\setminus \{ 0\} \right\} . 
\end{equation}
\end{lemma}

{\sl Proof: }
The proof uses classical arguments (see, e.g. \cite{FleRi}). 
We may restrict ourselves to nonnegative functions $u$.  
By (\ref{isopproblem}) and the coarea formula we obtain,
\begin{eqnarray}
\label{coarea1}
\int_{\mathbb{R} ^N } |x|^k |\nabla u| \, dx & = & 
\int _0 ^{\infty } \int\limits_{u=t } |x|^k \, \mathscr{H} _{N-1 } (dx) \, dt 
\\
\nonumber
 & \geq & C_{k,l,N} \int_0 ^{\infty } \left( \int_{u>t } |x|^l \, dx 
\right) ^{(k+N-1)/(l+N)} \, dt.
\end{eqnarray}
Further,
Cavalieri's principle gives
\begin{equation}
\label{cavalieri}
u(x)= \int_0 ^{\infty } \chi _{\{ u>t\} } (x)\, dt , \quad (x\in \mathbb{R} ^N ).
\end{equation}
Hence (\ref{cavalieri}) and Minkowski's inequality for integrals (see \cite{Stein}) lead to 
\begin{eqnarray}
\label{ineqmeas}   
 & & 
\\
\nonumber  
\int_{\mathbb{R} ^N } |x|^l |u|^{(l+N)/(k+N-1)} \, dx & = &
\int_{\mathbb{R} ^N } |x|^l \left| \int_0 ^{\infty} 
\chi_{\{ u>t\} } (x)\, dt \right| ^{(l+N)/(k+N-1)} \, dx 
\\
\nonumber 
 & \leq & \left( \int_0 ^{\infty } \left( 
\int_{\mathbb{R}^N } |x|^l \chi _{\{ u>t \} } (x) \, dx \right) 
^{(k+N-1)/(l+N)} \, dt \right) ^{(l+N)/(k+N-1)} 
\\
\nonumber
 & = & \left( \int_0 ^{\infty } \left( \int_{ u>t } |x|^l \, dx \right) 
^{(k+N-1)/(l+N)} dt \right) ^{(l+N)/(k+N-1)} .
\end{eqnarray}
Now (\ref{coarea1}) and (\ref{ineqmeas}) yield
\begin{equation}
\label{ineqQ1}
{\mathcal Q}_{k,l,N} (u) \geq C_{k,l,N} \quad \forall u\in C_0 ^1 (\mathbb{R} ^N ).
\end{equation}
To show (\ref{ineqQR}), 
let $\varepsilon >0$, and choose a smooth set 
$\Omega $ such that 
\begin{equation}
\label{ineqR1}
{\mathcal R}_{k,l,N} (\Omega ) \leq C_{k,l,N} +\varepsilon .
\end{equation}
It is well-known that there exists a sequence $\{ u_n \} \subset 
C_0 ^{\infty } (\mathbb{R} ^N )\setminus \{ 0\} $ 
such that 
\begin{eqnarray}
\label{limperim}
\lim_{n\to \infty } \int_{\mathbb{R} ^N } |x|^k |\nabla u_n | \, dx = 
\int_{\partial \Omega } |x|^k \, \mathscr{H} _{N-1} (dx) ,
\\
\label{limmeas}
\lim_{n\to \infty } \int_{\mathbb{R} ^N } |x|^l 
|u_n |^{(l+N)/(k+N-1)}  \, dx = \int_{ \Omega } |x|^l \, dx.
\end{eqnarray}
To do this, one may choose mollifiers of $\chi _{\Omega } $ 
as $u_n $ (see e.g. \cite{Talenti1}).  
Hence, for large enough $n$ we have
\begin{equation}
\label{ineqQ2}
{\mathcal Q}_{k,l,N} (u_n ) \leq C_{k,l,N} + 2\varepsilon .
\end{equation}
Since $\varepsilon $ was arbitrary, (\ref{ineqQR}) now 
follows from (\ref{ineqQ1}) and (\ref{ineqQ2}).      
$\hfill \Box $
%\color{red}

\begin{remark}
\rm 
Lemma \ref{R2} improves on \cite{AH}, Corollary 1.1 and 1.2, where the authors showed 
the inequality
$$
\inf \left\{ {\mathcal R}_{k,l,N} (\Omega ) :\,  \Omega \subset \mathbb{R}^N , \, 
\Omega \mbox{ smooth} \right\} \geq
\inf \left\{ {\mathcal Q}_{k,l,N} (u) :\, u\in C_0 ^1 (\mathbb{R} ^N )\setminus \{ 0\} , 
\right\} .
$$
\end{remark} 
%\color{black}

\section{Necessary conditions}

Throughout this section we assume that assumptions \eqref{ass1} are fullfilled, i.e.
\begin{equation*}
k+N-1 >0 \ \mbox{ and } \ l+N>0 .
\end{equation*}
The main result in this section is Theorem \ref{R4}  which concerns a phenomenon of symmetry breaking.

\noindent The following result holds true.

\begin{lemma}
\label{R3}
A necessary condition for 
\begin{equation}
\label{C>0}
C_{k,l,N} >0
\end{equation}
 is 
\begin{equation}
\label{k_l_ineq1}
 l \frac{N-1}{N} \leq k .
\end{equation}
\end{lemma}

{\sl Proof:} Assume that $k<l(N-1 )/N $, and let 
$te_1 = (t, 0, \ldots , 0)$, ($t>2$). It is easy to see that
there is a positive constant $D= D(k,l, N) $ such that
$$
{\mathcal R}_{k,l,N} (B_1 (te_1 ) ) \leq 
D \frac{ t ^k}{ t ^{l (k+N-1)/(l+N)} } ,
$$
 Since $k-l(k+N-1)/(l+N) <0$, it follows that
$$
\lim_{t\to \infty } {\mathcal R}_{k,l,N} (B_1 (te_1 ) ) =0.
$$
$\hfill \Box $

\begin{theorem}
\label{R4} 

A necessary condition for 
\begin{equation}
\label{isop1}
 C_{k,l,N} = C_{k,l,N} ^{rad} 
\end{equation}
is 
\begin{equation}
\label{k_l_ineq2}
l+1 \leq k + \frac{N-1}{ k+N-1} .
\end{equation}
\end{theorem}
\medskip

\begin{remark}\rm
Theorem \ref{R4} means that if  $l+1 > k + \frac{N-1}{ k+N-1} $,  then symmetry breaking occurs, that is $C_{k,l,N} < C_{k,l,N} ^{rad}$.
\end{remark}

{\sl Proof:} First we assume $N\geq 2$. Let $(r, \theta )$ denote 
$N$--dimensional spherical coordinates, 
$u \in C^2 (\mathscr{S}^{N-1} )$, $s\in C^2 (\mathbb{R})$ with $s(0)=0$, 
and define
$$
U(t ) := \{ x=r\theta \in \mathbb{R} ^N : \, 0\leq r < 1+ t u(\theta ) + s(t) \} , 
\quad (t\in \mathbb{R} ).
$$
Note that $U(0)= B_1 $. 
By the Implicit Function Theorem, we may choose $s$ in such a way that 
\begin{equation}
\label{intid1}
\int_{U(t)} |x|^l \, dx = \int _{B_1 } |x|^l \, dx \quad \mbox{for $|t|<t_0$},
\end{equation}
for some number $t_0 >0$. We set $s_1 := s'(0) $ and $s_2 := s^{\prime \prime} (0)$. 
Since
$$
\int_{U(t)} |x|^l \, dx = \int_{\mathscr{S} ^{N-1 } } 
\int_0 ^{1+ t u(\theta ) + s(t)}  \rho ^{l+N-1} \, d\rho \, d\theta,
$$
a differentiation of (\ref{intid1}) leads to
\begin{eqnarray}
\label{intid2}
0 & = & \int_{\mathscr{S} ^{N-1} } (u+ s_1 )\, d\theta \quad \mbox{and }
\\
\label{intid3}
0 & = & (l+N-1) \int_{\mathscr{S} ^{N-1} } (u+ s_1 )^2 \, d\theta + s_2 
\int_{\mathscr{S} ^{N-1} }  \, d\theta .
\end{eqnarray}
Next we consider the perimeter functional 
\begin{eqnarray}
\label{perim}
J(t) & := & \int_{\partial U(t)} |x|^k \, dx 
\\
\nonumber
 & = & \int_{\mathscr{S} ^{N-1} } (1+tu + s(t) )^{k+N-2} 
\sqrt{ (1+ tu+s(t) )^2 + t^2 |\nabla _{\theta } u|^2 } \, d\theta ,
\end{eqnarray}
where $\nabla _{\theta }$ denotes the gradient on the sphere. 
Differentiation of (\ref{perim}) leads to 
\begin{eqnarray*}
J'(0) & = & (k+N-1) \int_{\mathscr{S} ^{N-1} } (u+ s_1 )\, d\theta , \quad \mbox{and }
\\
J^{\prime \prime } (0) & = & 
(k+N-2) (k+N-1) \int_{\mathscr{S} ^{N-1} }  (u+s_1 )^2 \, d\theta + 
\\
 & & + (k+N-1) s_2 \int_{\mathscr{S} ^{N-1} } \, d\theta +  \int_{\mathscr{S} ^{N-1} } 
|\nabla _{\theta } u|^2 \, d\theta .
\end{eqnarray*}
By (\ref{intid2}) and (\ref{intid3}) this implies 
\begin{eqnarray}
\label{Jprime}
J'(0) & = & 0, \quad \mbox{and}
\\
\label{Jprimeprime}
J^{\prime \prime } (0)
 & = & 
(k+N-1) (k-l-1) \int_{\mathscr{S} ^{N-1} }  (u+s_1 )^2 \, d\theta + 
 \int_{\mathscr{S} ^{N-1} } |\nabla _{\theta } u|^2 \, d\theta .
\end{eqnarray} 
Now assume that (\ref{isop1}) holds. 
Then we have ${\mathcal R}_{k,l,N} (U(t)) \geq {\mathcal R}_{k,l,N} (B_1 )$ 
for all $t$ with $|t|<t_0 $. In view of (\ref{intid1}) 
this means that $J(t) \geq J(0) $ for $|t|<t_0 $, that is, 
\begin{equation}
\label{Jderiv}
J^{\prime \prime } (0) \geq 0 = J'(0).
\end{equation}
The second condition is (\ref{Jprime}), and the first condition 
implies in view of (\ref{intid2}) and 
(\ref{Jprimeprime}),
\begin{eqnarray}
\label{intineq1}
0 & \leq & (k+N-1)(k-l-1) \int_{\mathscr{S}^{N-1} }   v^2 \, d\theta + 
 \int_{\mathscr{S} ^{N-1} } |\nabla _{\theta } v|^2 \, d\theta 
\\
\nonumber
 & & \forall v\in C^2 (\mathscr{S} ^{N-1} ) \ \mbox{ with } \ 
 \int_{\mathscr{S} ^{N-1} } \, d\theta =0.
\end{eqnarray}
Let $V$ be the first non-trivial eigenfunction of the Laplace-Beltrami 
operator on the sphere. Then 
$
\int_{\mathscr{S}^{N-1 } } |\nabla _{\theta } V |^2 \, d\theta = 
(N-1) \int_{\mathscr{S} ^{N-1} } V^2 \, d\theta $ and 
$\int_{\mathscr{S} ^{N-1 }} V\, d\theta =0$. 
Choosing $v=V$ in (\ref{intineq1}), we obtain 
(\ref{k_l_ineq2}).
\\
Next assume that $N=1$. We proceed similarly as before. 
Let  $s\in C^2 (\mathbb{R} )$ with $s(0)=0$ and 
$U(t) := (-1+t , 1+ s(t) )$, ($t\in \mathbb{R}$). Note that $U(0)=(-1,1)= B_1 $.  
We may choose $s$ in such a way that 
\begin{equation}
\label{intid1N=1}
\mu _l (U(t)) =\mu_l (B_1 ) \quad \mbox{for $|t|< t_0 $}.
\end{equation} 
Setting $s_1 := s^{\prime }(0) $ and $s_2 := s^{\prime \prime } (0)$, 
a differentiation of (\ref{intid1N=1}) yields 
\begin{equation}
\label{intid2N=1}
s_1 =1 \ \mbox{ and } \ s_2 = -2l.
\end{equation}  
Next, let 
\begin{equation}
\label{perimN=1}
J(t) := P_{\mu _k } (U(t)) = |1+ s(t)| ^k + |-1+t|^k .
\end{equation}
A differentiation of this gives
\begin{equation}
\label{Jder}
J^{\prime } (0) =k(-1+s_1 ) =0 \ \mbox{ and } \ J^{\prime \prime } (0) 
= k (2k-2+s_2 ) = 2k(k-1-l).
\end{equation}
As before, we must have $ J^{\prime \prime } (0)\geq 0$,
so that (\ref{Jder}) implies $l+ 1 \leq k$.   
$\hfill \Box $

\section{Main results}

This section is devoted to the proof of Thorem \ref{maintheorem}, that is we obtain sufficient conditions on $k,l$ and $N$ such that 
$ C_{k,l,N} = C_{k,l,N} ^{rad}$ holds, or equivalently,
\begin{equation}
\label{ineqrad}
{\mathcal R}_{k,l,N} (M) \geq C_{k,l,N}^{rad}  
\quad \mbox{for all measurable sets $M$ with $0< \mu _l (M) <+\infty $.}
\end{equation}
Such a proof is contained in various subsections each of one addresses one of the cases of Theorem 1.1.

\noindent Throughout this section  we again assume (\ref{ass1}), i.e.
\begin{equation*}
k+N-1 >0 \ \mbox{ and } \ l+N>0 .
\end{equation*}

\medskip

\subsection{Proof of Theorem \ref{maintheorem}, case (i).}
As mentioned in the Introduction, Theorem \ref{maintheorem} was already shown under assumption {\bf (i)} in \cite{Howe}. 
Below we give another simple proof which is based on Gauss' Divergence Theorem. Note that this tool has been applied in similar situations in 
\cite{KZ} and \cite{BrasPhil}. We also discuss equality cases of (\ref{mainineq}).  
\begin{theorem}    
\label{R5}
Let $l+1 \leq k $. Then (\ref{isop1}) holds. 
Moreover, if $l+1 <k$ and
\begin{equation}
\label{M=BR}
{\mathcal R}_{k,l,N} (M) = C_{k,l,N} ^{rad} \ 
\mbox{ for some measurable set $M$ with $0<\mu _l (M)< +\infty $},
\end{equation}  
then $M = B_R $ for some $R>0$. 
\end{theorem}
{\sl Proof: } We consider two cases. 
\\
{\bf 1.} $l+1=k$.  
\\
Let $\Omega $ be smooth. We choose $R>0$ such that $\Omega ^{\star} = B_R $.
From Gauss' Divergence Theorem we have, 
($\nu $ denotes the exterior unit normal to $\partial \Omega $),
\begin{eqnarray}
\label{gauss}
\int _{\Omega } |x|^{l} \, dx
 & = & 
\frac{1}{l+N} \int_{\Omega } \mbox{div}\, \left( x|x|^{l} \right) \, dx 
\\
 \nonumber
 & = & 
\frac{1}{l+N} \int_{\partial \Omega } (x\cdot \nu )  |x|^{l} \mathscr{H}_{N-1}   ( dx) 
\\
 \nonumber
 & \leq & 
\frac{1}{l+N} \int_{\partial \Omega }  |x|^{l+1} \mathscr{H}_{N-1}   ( dx),
\end{eqnarray}
with equality for $\Omega = B_R $,
and (\ref{ineqrad}) follows for smooth sets. 
Using (\ref{CklNsmooth}), we also obtain (\ref{ineqrad}) for measurable sets. 
\\
{\bf 2.} $l+1<k$. 
\\ 
Using Lemma \ref{rangekl1} and the result for $l+1=k$ 
we again obtain (\ref{ineqrad}), 
and (\ref{M=BR}) can hold only if $M= B_R $ for some $R>0$.
$\hfill \Box $ 
\begin{corollary}
\label{sufficiency}
Condition (\ref{k_l_ineq1}), i.e. $l\frac{N-1}{N}\le k$ is a 
necessary and sufficient condition for 
$C_{k,l,N} >0$.
\end{corollary}
{\sl Proof:\/} 
The necessity follows from Lemma \ref{R3}, and the sufficiency in the case 
$l+1\leq k$ follows from Theorem \ref{R5}.
Finally, assume that $k< l+1$. Then  (\ref{isopproblem}) is equivalent to 
(\ref{ineqQR}), by Lemma \ref{R2}.
Now the main Theorem of \cite{CKN} tells us that condition 
(\ref{k_l_ineq1}) is also sufficient  
for $C_{k,l,N} >0$. 
$\hfill \Box $
\\[0.1cm]
\medskip

\subsection{Proof of Theorem \ref{maintheorem}, case (ii).} 
The case when $k$ assumes negative values has been settled in a recent paper, see \cite{ChiHo}, Theorem 1.3. We slightly improve on this result by adding a full treatment of the equality case in (\ref{isop1}). For the convenience of the reader, we include the full proof. 
\begin{theorem}
\label{th1bis} 
Let $N\geq 2 $, and let $k,l$ satisfy 
\begin{equation} \label{lk} 
l \frac{N-1}{N} \leq k 
\leq \min\{0, l+1\}. 
\end{equation}
Then (\ref{isop1}) holds. 
Moreover
(\ref{M=BR}) holds only if $M= B_R $ for some $R>0$. 
\end{theorem}
{\sl Proof  :\/} 
Let $u\in C^{\infty }_0(\mathbb{R} ^N)\setminus \{ 0 \} $.
We set
$$
y:=x|x|^\frac{k}{N-1}\, , \quad v(y):=u(x)\, , \quad 
s:=r^\frac{k+N-1}{N-1}\,
.
$$
Using $N$-dimensional spherical coordinates, let $\nabla_\theta$ denote  
the tangential part of the gradient on  
 ${\mathcal S^{N-1}}$. Then we obtain
\begin{eqnarray} 
\label{cambio1} 
 & &
\int_{\mathbb{R} ^N}
|x|^l
|u|^{(l+N)/(k+N-1)} \, dx
\\
\nonumber
  & = & 
\int_{\mathscr{S}^{N-1}} 
\int_0^{\infty} 
r^{l+N-1} |u|^{(l+N)/(k+N-1) }\, 
dr\, d\theta 
\\
\nonumber
 & = & 
 \frac{N-1}{k+N-1} 
 \int_{\mathscr{S}^{N-1}} 
 \int_0^{\infty} 
 s^{\frac{l+N}{k+N-1}(N-1)-1 }
 |v|^{(l+N)/(k+N-1)}\, 
 ds \, d\theta 
 \\
\nonumber
 & = & 
 \frac{N-1}{k+N-1} 
 \int_{\mathbb{R} ^N} 
 |y|^{\frac{l+N}{k+N-1}(N-1)-N}|v|^{(l+N)/(k+N-1)}\, dy  
 \\
\nonumber
 & = & 
 \frac{N-1}{k+N-1} 
 \int_{\mathbb{R} ^N} 
 |y|^{(l(N-1)-kN)/(k+N-1)}
 |v|^{(l+N)/(k+N-1)}\, dy
 \, .
\end{eqnarray}
Further we calculate
\begin{eqnarray}
\label{cambio2}
 \int_{\mathbb{R} ^N}|x|^k |\nabla_x u| \, dx
 & = &
 \int_{\mathscr{S}^{N-1} } 
 \int_0^{\infty} 
 r^{k+N-1}
 \left(  
 u_r ^2 +\frac{|\nabla_\theta  u|^2}{r^2}   
 \right) ^{1/2}\, 
 dr \, d\theta
 \\
 \nonumber
 & = &
 \int_{\mathscr{S}^{N-1} } 
 \int_0^{\infty} 
 s^{N-1}
 \left(  
 v_s ^2+\frac{|\nabla_\theta v|^2}{s^2} 
 \left(
 \frac{N-1}{k+N-1} \right) ^2  \right) ^{1/2} \, 
 ds \, d\theta
\nonumber 
\\
\nonumber
 & \geq &
 \int_{\mathscr{S}^{N-1} } 
 \int_0^{\infty} 
 s^{N-1}
 \left(  
 v_s ^2 +\frac{|\nabla_\theta v|^2}{s^2}   
 \right) ^{1/2} \, 
 ds \, d\theta
\\
\nonumber 
 & = & 
 \int_{\mathbb{R} ^N} |\nabla_y v| \, dy \, ,
\end{eqnarray}
where we have used (\ref{lk}).
By \eqref{cambio1} and \eqref{cambio2} we deduce,
 \begin{eqnarray}
 \label{Q2}
  & &
 {\mathcal Q}_{k,l,N}(u)
 \\
 \nonumber
  & \geq & 
 \frac{\displaystyle 
 \int_{\R ^N} |\nabla_y v| \, dy}{\displaystyle 
 \left( 
 \int_{\R ^N}|y|^{l' }|v|^{(l+N)/(k+N-1)}\, dy  \right) 
 ^{(k+N-1)/(l+N)} 
 }
 \left(
\frac{k+N-1}{N-1}
\right) ^{(k+N-1)/(l+N)}
\\
\nonumber
 & = &
\left( \frac{k+N-1}{N-1} \right)
^{(k+N-1)/(l+N)} 
{\mathcal Q}_{0,l' ,N }(v)\, ,
\end{eqnarray}
where we have set $l' :=\frac{l(N-1)-kN}{k+N-1}$. 
Note that we have $-1 \leq l' \leq 0$ by the assumptions (\ref{lk}).
\\
Hence we may apply Lemma \ref{R2} to both sides of (\ref{Q2}). 
This yields
\begin{equation}
\label{relationCC}
C_{k,l,N} \geq \left( \frac{k+N-1}{N-1} \right) ^{(k+N-1)/(l+N)} C_{0,l', N} .
\end{equation}
Furthermore, Lemma \ref{rangekl1} tells us that
\begin{equation}
\label{CCrad}
C_{0,l',N} = C_{0,l', N} ^{rad} .
\end{equation}
Since also 
$$
\left( \frac{k+N-1}{N-1} \right)
^{(k+N-1)/(l+N)} C_{0,l',N} ^{rad} =C_{k,l,N}^{rad} \, .
$$
From this, (\ref{relationCC}) and (\ref{CCrad}), 
we deduce that  $C_{k,l,N}\ge C_{k,l,N}^{rad}$. 
Since $C_{k,l,N}\le C_{k,l,N}^{rad}$ by definition, (\ref{isop1}) follows.
\\
Next assume that (\ref{M=BR}) holds.  
If $l(N-1)/N <k$, then Lemma \ref{rangekl1} tells us that we must have 
$M=B_R $ for some $R>0$.
Hence it remains to consider the case 
$$
l\frac{N-1}{N} =k<0.
$$
Then 
$$
l' =
\frac{l(N-1)-kN}{k+N-1} =0.
$$
Setting $k_1 := l(N-1)/N$ and 
$$
\widehat{M} := \{ y= x|x|^{\frac{l}{N} } :\, x\in M\} ,
$$
(\ref{Q2}), the classical isoperimetric inequality (\ref{isopclass})  
and a limit argument analogous to the proof of Lemma \ref{R2} leads to
\begin{eqnarray}
\nonumber
C_{k_1 ,l,N} ^{rad}
 & = & 
{\mathcal R}_{k_1 ,l,N}  (M)
\geq 
\left( \frac{l+N}{N} \right) ^{(N-1)/N} {\mathcal R}_{0,0,N} (\widehat{M})
\\
\label{compareCC1}
 & \geq  &
 \left( \frac{l+N}{N} \right) ^{(N-1)/N} C_{0,0,N} ^{rad} .
\end{eqnarray}
Since 
$$
C_{k_1 ,l,N} ^{rad} = \left( \frac{l+N}{N} \right) ^{(N-1)/N} C_{0,0,N} ^{rad} ,
$$ 
(\ref{compareCC1}) implies that
${\mathcal R}_{0,0,N} (\widehat{M}) = C_{0,0,N} ^{rad} $. 
By (\ref{isopclass}) it follows that $\widehat{M} = B_{\widehat R} (y_0 ) $ 
for some $\widehat{R} >0 $ and $y_0 \in \mathbb{R}^N $. 
Using again (\ref{compareCC1}) we find  
$$
{\mathcal R}_{k_1 ,l,N} (M) = \left( \frac{l+N}{N} \right) ^{(N-1)/N} 
{\mathcal R}_{0,0,N} (B_{\widehat{R}} (y_0 )).
$$
Since 
$$ 
\left( \frac{l+N}{N} \right) ^{(N-1)/N} \mu _l (M) = 
\mu _0 (B_{\widehat{R} } (y_0 ) ),
$$ 
this implies
$$
P_{\mu _{k_1 } } (M) = P_{\mu _0 } (B_{\widehat{R} } (y_0 )).
$$
It is easy to see that this is possible only when 
$y_0 =0$.
$\hfill \Box $

\begin{remark} 
\rm

{\bf (a)} A well-known special case of Theorem \ref{th1bis} is $k=0 \geq l >-N $, 
see for instance \cite{M}, p.11.
\\ 
{\bf (b)}
The idea to use spherical coordinates, 
and in particular the inequality (\ref{cambio2}) in our last proof, 
appeared already in some work of T. Horiuchi, 
see \cite{H} and \cite{HK}. 
\end{remark} 
\medskip

\subsection{Proof of Theorem \ref{maintheorem}, case (iii).} 
Now we treat the case when $k$ assumes non-negative values. 
Throughout this subsection we assume  $N\geq 2$ and $k\leq l+1$. 
The main result is Theorem \ref{th1ter}. 
Its proof is long and  requires some auxiliary results. 
But the crucial idea is an interpolation argument that occurs in the proof of the following Lemma \ref{4.3}, formula (\ref{ineq2}). 
\begin{lemma}  
\label{4.3}
Assume $l(N-1)/N\leq k$ and $k\geq 0$.   
Let $u\in C_0 ^1 (\mathbb{R}^N)\setminus \{ 0 \} $, $u\geq 0$, 
and define $v$ by   
\begin{equation}
\label{transf1}
 v(y) := 
  u(x) \quad \mbox{ where 
$ y  :=  x |x|
 ^{ \frac{k}{N-1} }$,  ($x\in \mathbb{R} ^N $).}
\end{equation}
Then for every 
$A\in \left[ 
0, \frac{(N-1) ^2}{ (k+N-1 )^2 } 
\right] 
$,
\begin{equation}
\label{ineq1}
 {\mathcal Q}_{k,l,N} (u)  \geq \left( 
\frac{k+N-1}{N-1} 
\right) 
^{ \frac{k+N-1}{l+N} } 
 \cdot
\frac{ 
\left( 
\dint_{\mathbb{R} ^N } 
|\nabla _y v| \, 
\, dy 
\right) 
^A
\cdot 
\left( 
\dint_{\mathbb{R} ^N } 
| v_z|  
\, dy  
\right) 
^{1-A}
 }{ 
\left( 
\dint_{\mathbb{R} ^N } 
|y|
^{ \frac{l(N-1)-kN}{k+N-1} } 
v 
^{ \frac{l+N}{k+N-1} } 
\, dy 
\right) 
^{ \frac{k+N-1}{l+N} }
} 
.
\end{equation}
\end{lemma}    

{\sl Proof:} Setting 
$$
z:= |y| = r^{\frac{k+N-1}{N-1} },
$$ 
we calculate as in the proof of Theorem \ref{th1bis} ,
$$
\int_{\mathbb{R} ^N } |x| ^k |\nabla _x u | \, dx 
 =    
 \int_{{\mathscr S} ^{N-1} } 
 \int_0^{+\infty } z^{N-1}  
 \sqrt{ v_z ^2 + \frac{|\nabla _{\theta } v|^2 }{z^2 } 
 \frac{(N-1)^2}{(k+N-1)^2 } } 
 \, dz\, d\theta .
 $$
 Since the mapping
 $$
 t\longmapsto \log 
 \left(
 \int _{{\mathscr S} ^{N-1} } 
 \int_0^{+\infty } z^{N-1}  
 \sqrt{ v_z ^2 + t\frac{|\nabla _{\theta } v|^2 }{z^2 }  } \, dz\, d\theta \right)
 $$
 is concave, we deduce that for every 
 $A\in \left[ 0, \frac{(N-1) ^2}{ (k+N-1 )^2 } \right] $,
\begin{eqnarray}
\label{ineq2}
 & & \int_{\mathbb{R} ^N } |x| ^k |\nabla _x u | \, dx 
\\ 
\nonumber
 & \geq & 
 \left( 
 \int_{\mathscr{S}^{N-1} } \int_0 ^{+\infty } z^{N-1}  
 \sqrt{ 
 v_z ^2 + \frac{|\nabla _{\theta } v|^2 }{z^2 } 
  } 
  \, dz\, d\theta \right)
  ^A 
  \cdot 
\left(  
\int_{\mathscr{S}^{N-1}}
\int_0^{+\infty } z^{N-1}  |v_z |  \, dz\, d\theta 
\right) 
^{1-A}
\\
\nonumber
 & = & 
 \left( 
 \int_{\mathbb{R}^N } |\nabla _y v| \, dy 
 \right) 
 ^A
  \cdot 
 \left( 
 \int_{\mathbb{R}^N } |v_z | \, dy 
 \right) 
 ^{1-A} .
 \end{eqnarray}
 Finally, we have 
 \begin{equation}
 \label{equaldenom}
 \int_{\mathbb{R}^N } |x| ^l u
 ^{ \frac{l+N}{k+N-1} } \, dx   = 
 \frac{N-1}{k+N-1} \int_{\mathbb{R}^N} |y| 
 ^{ \frac{l(N-1)-kN}{k+N-1} } 
 v
 ^{ \frac{l+N}{k+N-1} } \, dy .
 \end{equation}
Now (\ref{ineq1}) follows from 
(\ref{ineq2}) and (\ref{equaldenom}).
$\hfill \Box $
\\[0.1cm]   
Next we want to estimate the right-hand-side of (\ref{ineq1}) from below. 
We will need a few more properties of the starshaped rearrangement.
\begin{lemma}  
\label{4.2}  
Assume $l(N-1)/N\leq k$. 
Then we have for any function 
$v\in C_0 ^1 (\mathbb{R}^N )\setminus \{ 0 \}$ with $v\geq 0$, 
\begin{eqnarray}
\label{starsh4}
 & & \int_{\mathbb{R}^N }  v
 ^{ \frac{N}{N-1} } 
\, dy 
= 
\int_{\mathbb{R}^N }  \widetilde{v}
^{ \frac{N}{N-1} } 
\, dy 
\\
\label{starsh5}
 & & 
 \int_{\mathbb{R}^N } |y| 
 ^{ \frac{l(N-1)-kN}{k+N-1} } 
 v
 ^{ \frac{l+N}{k+N-1} } 
 \, dy 
\leq 
\int_{\mathbb{R}^N } |y| 
^{ \frac{l(N-1)-kN}{k+N-1} } 
\widetilde{v}
^{ \frac{l+N}{k+N-1} } 
\, dy,
\\
\label{vL1}
 & & 
\frac{ 
y\cdot \nabla \widetilde{v} }{|y|} 
\equiv 
\frac{
\partial \widetilde{v} }{
\partial z } 
\in L^1 (\mathbb{R} ^N )
\quad \mbox{and }
\\
\label{starsh6}
 & &
\int_{\mathbb{R}^N } 
\left| \frac{ \partial v}{\partial z}  \right| \, dy 
\geq  
\int_{\mathbb{R}^N } 
\left| \frac{ \partial \widetilde{v} }{\partial z} \right|
\, dy.
\end{eqnarray}
\end{lemma}

{\sl Proof:} 
Equality (\ref{starsh4}) follows from (\ref{caval1}). 
Now let us prove (\ref{starsh5}). Set 
$$
w(y):= |y| 
^{ \frac{l(N-1)-kN}{l+N} } . 
$$
Since $l(N-1)-kN \leq 0$, we have $w= \widetilde{w} $. 
Hence (\ref{starsh5}) follows from (\ref{harlit}) and (\ref{monrearr}). 
\\
Next let $\zeta := z^N  $ and define $V$ and $\hat{V} $ by 
$V(\zeta ,\theta) := v(z\theta )$, and 
$\widehat{V} (\zeta ,\theta ) := \widetilde{v} (z\theta )$. 
Observe that for each $\theta \in \mathscr{S}^{N-1} $,  
$\widehat{V} (\cdot , \theta ) $ is the equimeasurable 
non-increasing rearrangement of 
$V (\cdot ,\theta )$. Further we have 
$$
\frac{ \partial v }{ \partial z }
  =  
 N\zeta 
 ^{ \frac{N-1}{N}  } 
 \frac{ \partial V}{\partial \zeta } \ \mbox{ and } \ 
\frac{ \partial \widetilde{v} }{ \partial z }
  = 
 N\zeta 
 ^{ \frac{N-1}{N}  } 
 \frac{ \partial \widehat{V} }{ \partial \zeta  } 
 .
$$
Since $\frac{\partial v}{\partial z} \in L^{\infty } (\mathbb{R}^N )$, 
Lemma \ref{Landes} tells us 
that for every $\theta \in {\mathscr S} ^{N-1} $,
\begin{eqnarray*} 
\int_0^{+\infty } z^{N-1} \left| \frac{\partial v}{\partial z} 
(z\theta ) \right| \, dz 
 & = & 
 \int_0 ^{+\infty } \zeta 
 ^{ \frac{N-1}{N} } 
 \left| \frac{\partial V}{\partial \zeta } (\zeta ,\theta )
 \right| \, d \zeta  
 \\
 & \geq & 
 \int_0 ^{+\infty } \zeta 
 ^{ \frac{N-1}{N} } 
 \left| \frac{\partial \widehat{V} }{\partial \zeta }  
 (\zeta ,\theta )\right| \, d\zeta 
 \\
  & = & 
\int_0^{+\infty } z^{N-1} \left|
\frac{ \partial \widetilde{v} }{\partial z} (z\theta )\right| \, dz .
\end{eqnarray*}
Integrating this over ${\mathscr S} ^{N-1}$, 
we obtain (\ref{starsh6}).
$\hfill \Box$

A final ingredient is
\begin{lemma}  
\label{4.1}
Assume that $l(N-1)/N\leq k$, 
and let $M$ be a bounded starshaped set. Then
\begin{eqnarray}
\label{holder1}
 & & 
 \left( 
 \int_M |y| 
 ^{ \frac{l(N-1) -kN}{k+N-1} } 
 \, dy 
 \right) 
 ^{ \frac{k+N-1}{l+N} } 
 \\
 \nonumber
  & \leq & 
  d_1 
  \left( 
  \int_M \, dy 
  \right) 
  ^{  \frac{(N-1)(l-k+1) }{ l+N}   } 
  \cdot 
  \left( 
  \int_M |y|^{-1} \, dy 
  \right) 
  ^{  \frac{kN-l(N-1) }{ l+N}  } ,
\quad \mbox{ where}
\\
 & & 
\label{d1}
d_1 = 
\left( 
\frac{k+N-1}{l+N} 
\right) 
^{ \frac{k+N-1}{l+N} } 
\cdot 
\left( 
\frac{N}{N-1} 
\right) 
^{ \frac{(N-1)(l-k+1)}{l+N} } .
\end{eqnarray}
Moreover,
if $k<l+1$ and $l(N-1)/N <k$, then 
equality in (\ref{holder1}) holds only
if $M=B_R $ for some $R>0$.
\end{lemma}    

{\sl Proof:} Since $M$ is starshaped, 
there is a bounded measurable function 
$m : \mathscr{S} ^{N-1} \to [0, +\infty )$, such that
\begin{equation} 
\label{Mrepresent}
M= \{ z\theta :\, 0\leq z < m (\theta ), \ 
\theta \in \mathscr{S} ^{N-1} \} .
\end{equation} 
Using H\"older's inequality we obtain
\begin{eqnarray}
\label{chain}
 & & \int_M |y| 
 ^{ \frac{l(N-1) -kN}{k+N-1} } 
 \, dy 
 \\
 \nonumber
 & = & 
 \frac{k+N-1}{(l+N)(N-1)} 
\int_{{\mathscr S} ^{N-1} } m (\theta ) 
^{ \frac{(l+N)(N-1)}{k+N-1} } 
\, d\theta 
\\
 \nonumber
 & = & 
 \frac{k+N-1}{(l+N)(N-1)} 
\int_{{\mathscr S} ^{N-1} } m (\theta ) 
^{ \frac{kN-l(N-1)}{k+N-1}(N-1) } m (\theta ) 
^{ \frac{(N-1)(l-k+1)}{k+N-1} N} 
\, d\theta 
\\
\nonumber
& \leq &    
 \frac{k+N-1}{(l+N)(N-1)} 
 \left( 
\int_{{\mathscr S} ^{N-1} } 
m (\theta ) ^N \, d \theta 
\right) 
^{ \frac{(N-1)(l-k+1)}{k+N-1} } 
\cdot 
\left( 
\int_{{\mathscr S} ^{N-1} } 
m (\theta ) ^{N-1} \, d \theta 
\right) 
^{ \frac{kN- l(N-1)}{k+N-1} }
\\
\nonumber
 & = &  
 \frac{k+N-1}{(l+N)(N-1)} 
 \left( 
 N \int_M dy 
 \right) 
 ^{ \frac{(N-1)(l-k+1)}{k+N-1} } 
 \cdot 
 \left( 
 (N-1) \int_M |y| ^{-1} \, dy 
 \right) 
 ^{ \frac{kN- l(N-1)}{k+N-1} } ,
\end{eqnarray}
and (\ref{holder1}) follows.
If $k<l+1$ and $l(N-1)/N< k$, then (\ref{chain}) holds   
with equality only if $m (\theta )=\mbox{const }$.   
$\hfill \Box $
\\[0.1cm] 
 Now we are ready to prove our main result. 
\begin{theorem}   
\label{th1ter}
Assume $N\geq 3$, $0\le k\leq l+1$ and
\begin{equation}
\label{crucial}
\frac{1}{l+N} \geq 
\frac{1}{k+N-1} - \frac{(N-1)^2 }{ N(k+N-1)^3 } .
\end{equation}
Then (\ref{isop1}) holds.
Furthermore, if 
 inequality (\ref{crucial}) is strict, or if $k>0$,
then (\ref{M=BR}) holds only if $M=B_R $ for some $R>0$.  
\end{theorem}    

{\sl Proof: } First observe that the conditions 
$k\geq 0$ and (\ref{crucial}) also imply
$l(N-1)/N \leq k$. Let $u \in C_0 ^{\infty } (\mathbb{R}^N)\setminus \{ 0\} $, $u\geq 0$,  
and let $v$ be given by 
(\ref{transf1}).
In view of (\ref{crucial}), we may choose
$$
A=\frac{N(l-k+1)}{l+N} 
$$
to obtain 
\begin{equation}
\label{ineq5bis}
 {\mathcal Q}_{k,l,N} (u)
\geq 
 \left( 
 \frac{k+N-1}{N-1} 
 \right) 
^{ \frac{k+N-1 }{ l+N} } 
\cdot 
\frac{ 
\left( 
\dint_{\mathbb{R} ^N } 
|\nabla _y v|  
\, dy 
\right) 
^{ \frac{N(l-k+1) }{ l+N} } 
\cdot 
\left( 
\dint_{\mathbb{R} ^N } 
|v_z |
\, dy  
\right) 
^{ \frac{kN-l(N-1) }{ l+N} } 
 }{
\left( 
\dint_{\mathbb{R} ^N } |y|
^{ \frac{l(N-1)-kN }{ k+N-1} }
v 
^{ \frac{l+N}{k+N-1} } 
\, dy 
\right) 
^{ \frac{k+N-1}{l+N} } 
 }  .
\end{equation} 
Further, (\ref{starsh6}) and 
Hardy's inequality yield
\begin{equation}
\label{ineq4}
 \int_{\mathbb{R}^N } |v_z| \, dy
  \geq   
  \int_{\mathbb{R}^N } |\widetilde{v}_z| \, dy
 \geq  (N-1) 
  \int_{\mathbb{R}^N } \frac{\widetilde{v}}{|y|}  
  \, dy.
 \end{equation}
Together with (\ref{ineq5bis}) and (\ref{starsh5}) this leads to
\begin{eqnarray}
\label{ineq5final}
 {\mathcal Q}_{k,l,N} (u) 
 & \geq & 
(N-1)
^{ \frac{kN-l(N-1)}{l+N} } 
 \left( 
 \frac{k+N-1}{N-1} 
 \right) 
^{ \frac{k+N-1 }{ l+N} } 
\cdot 
\\
\nonumber
 & & \cdot
\frac{ 
\left( 
\dint_{\mathbb{R} ^N } 
|\nabla _y v|  
\, dy 
\right) 
^{ \frac{N(l-k+1) }{ l+N} } 
\cdot 
\left( 
\dint_{\mathbb{R} ^N } 
\frac{\widetilde{v} }{|y|}  
\, dy  
\right) 
^{ \frac{kN-l(N-1) }{ l+N} } 
 }{
\left( 
\dint_{\mathbb{R} ^N } |y|
^{ \frac{l(N-1)-kN }{ k+N-1} }
\widetilde{v} 
^{ \frac{l+N}{k+N-1} } 
\, dy 
\right) 
^{ \frac{k+N-1}{l+N} } 
 }  .
\end{eqnarray} 
Now let $M $ be a bounded measurable set.
Then combining (\ref{limperim}), (\ref{limmeas}) 
and the argument leading to (\ref{CklNsmooth}) we deduce that  
there exists a sequence of non-negative functions 
$\{ u_n \} \subset C_0 ^1 (\mathbb{R}^N  )$ such that
\begin{equation}
\label{lim1}
\lim_{n\to \infty } 
\int_{\mathbb{R} ^N } |x| ^k |\nabla u_n | \, dx = 
P_{\mu _k } (M) 
\end{equation}
and
\begin{equation}
\label{lim2}
 u_n \longrightarrow \chi_{M } \quad 
\mbox{ in $L^p (\mathbb{R}^N ) $ for every $p\geq 1 $.}
\end{equation}
We define  
$M ':= \{ y= x|x|^{\frac{k}{N-1}} : \, x \in M \} $ 
and $v_n (y) := u_n (x) $, 
($y= x|x|^{\frac{k}{N-1}}$, $x\in \mathbb{R}^N  $).
Let $\widetilde{v_n } $ and $\widetilde{M '} $ 
be the starshaped rearrangements of $v_n $ and $M ' $ respectively. 
Then (\ref{lim1}) and (\ref{lim2}) also imply 
\begin{eqnarray}
\label{lim3}
 & & \lim_{n\to \infty } \int_{\mathbb{R}^N } |\nabla _y v_n | \, dy  = 
P_{\mu _0 } (M'  ),
\quad \mbox{and}
\\ 
\label{lim4}
 & & \widetilde{v_n }\longrightarrow 
\chi _{\widetilde{M ' } } \
\mbox{ in $L^p (\mathbb{R}^N ) $ for every $p\geq 1 $.}
\end{eqnarray}
Choosing $u=u_n $ in (\ref{ineq5final}) and passing to the limit 
$n\to \infty $, 
we obtain, using (\ref{lim1}), (\ref{lim2}), (\ref{lim3}), 
(\ref{lim4}) and the isoperimetric inequality (\ref{isopclass}),    
\begin{eqnarray}
\label{ineq6}
  {\mathcal R}_{k,l,N} (M )& \geq &  
(N-1)
^{ \frac{kN-l(N-1)}{l+N} } 
\left( 
 \frac{k+N-1}{N-1} 
\right) 
^{ \frac{k+N-1}{N-1} } 
\cdot 
\\
\nonumber
 & & \cdot
\frac{ 
\left( P_{\mu _0 } (M') \right)
^{ \frac{N(l-k+1)}{l+N} } 
\cdot 
\left( 
\dint_{\widetilde{M' } }  
\frac{dy}{|y| } 
\right) 
^{ \frac{kN-l(N-1)}{l+N} } 
}{
\left( 
\dint_{\widetilde{M '}  } |y|
^{ \frac{l(N-1)-kN}{k+N-1} } 
\, dy 
\right) 
^{ \frac{k+N-1}{l+N} } 
} 
\\
 \nonumber
 & \geq &  
(N-1)
^{ \frac{kN-l(N-1)}{l+N} } 
\left( 
N\omega_N ^{1/N} 
\right) 
^{ \frac{N(l-k+1)}{l+N} } 
\left( 
 \frac{k+N-1}{N-1} 
\right) 
^{ \frac{k+N-1}{N-1} } 
\cdot 
\\
\nonumber
 & & \cdot
\frac{ 
\left( \mu_0 (M') \right)
^{ \frac{(N-1)(l-k+1)}{l+N} } 
\cdot 
\left( 
\dint_{\widetilde{M' } }  
\frac{dy}{|y| } 
\right) 
^{ \frac{kN-l(N-1)}{l+N} } 
}{
\left( 
\dint_{\widetilde{M '}  } |y|
^{ \frac{l(N-1)-kN}{k+N-1} } 
\, dy 
\right) 
^{ \frac{k+N-1}{l+N} } 
} .
\end{eqnarray}
In view of (\ref{holder1}) and since 
$\mu _0 (M') = \mu_0 (\widetilde{M'})$
we finally get from this
\begin{eqnarray}
\label{ineq7}
{\mathcal R}_{k,l,N} (M) 
 & \geq &  
(N-1)
^{ \frac{kN-l(N-1)}{l+N} } 
\left( 
N\omega_N ^{1/N} 
\right) 
^{ \frac{N(l-k+1)}{l+N} }  
\left( \frac{k+N-1}{N-1} \right) 
^{ \frac{k+N-1}{N-1} } \frac{1}{d_1} 
\\
\nonumber 
 & = &  \left( N\omega _N \right) ^{\frac{l-k+1}{l+N} } 
\cdot (l+N)^{\frac{k+N-1}{l+N} } = C_{k,l,N} ^{rad}  ,
\end{eqnarray}
and (\ref{isop1}) follows by (\ref{CklNsmooth}).
\\
Now assume that (\ref{M=BR}) holds. 
If inequality (\ref{crucial}) is strict, 
then Lemma \ref{rangekl1} tells us that we must have 
$M= B_R $ for some $R>0$.
It remains to consider the case that $k>0$ and 
$$
\frac{1}{l+N} = 
\frac{1}{k+N-1} - \frac{(N-1)^2 }{ N(k+N-1)^3 } .
$$
Then we also have $l(N-1)/N <k $ and $k<l+1$. 
Now observe that all inequalities in (\ref{ineq6}) and (\ref{ineq7}) 
become equalities.
First, combining (\ref{isopclass}) and (\ref{ineq6}), we obtain that 
$M'= B_R (x_0 )$ 
for some $R>0 $ and   $x_0 \in \mathbb{R}^N $. Further, 
(\ref{holder1}) together with (\ref{ineq7}) imply that
$\widetilde{M'} $ is a ball centered at the origin. 
But this is possible only if $x_0 =0$.     
$\hfill \Box $
\begin{remark}  
\rm 
Theorem \ref{th1ter} is valid for $N\geq 2$.  Moreover, when 
$N\geq 3$, then (\ref{crucial}) covers the important range 
$$
l=0\leq k\leq 1. 
$$
However, we emphasize that this is not true in the case $N=2$  
(see however Lemma \ref{csato} in the next subsection). 
\end{remark}    
\medskip
\subsection{Proof of Theorem \ref{maintheorem}, case (iv).}
Next we improve on the subsections {\bf 5.2} and {\bf 5.3} in the two-dimensional case. 
We will make use of the following result of 
G. Csat\'o \cite{C}, that has been obtained by using conformal mappings.
\begin{lemma}
\label{csato}
 Let $N=2$, $l=0$ and $0\leq k\leq 1 $. Then (\ref{isop1}) holds. 
\end{lemma}  

The following result holds
\begin{corollary}
\label{N=2_l<0}
Let $N=2$, $k\leq l+1$, 
\begin{eqnarray} 
\label{lkineq}
 & &  \frac{l}{2} \leq k \quad \mbox{and}
\\
\label{lineq}
 & &  l\leq 0 .
\end{eqnarray}
Then (\ref{isop1}) holds. 
Furthermore, if $\frac{l}{2} <k$, then equality in (\ref{ineqrad}) 
holds only if $M=B_R $ for some $R>0$.  
\end{corollary}

{\sl Proof:} If $k\leq 0$, then (\ref{isop1}) follows from Theorem \ref{th1bis}. 
If $k\geq 0$, then (\ref{isop1}) 
follows from Lemma \ref{csato} together with Lemma \ref{rangekl1}.  
\\
Finally, assume that (\ref{M=BR}) holds and that $\frac{l}{2} <k$. 
Then the above result for $\frac{l}{2} =k$ and Lemma \ref{rangekl1} 
shows that $M= B_R $ for some $R>0 $.   
$\hfill \Box $
\\[0.1cm]
\begin{theorem}
\label{theoN=2k>0}
Let $N=2$, $k\leq l+1$,  
\begin{eqnarray}
\label{k>0N2}
 & & 
k\geq \frac{1}{3} \quad \mbox{and}
\\
\label{crucial2}
 & & 
\frac{1}{l+2} \geq \frac{1}{k+1} -\frac{16}{27 (k+1)^3 } .
\end{eqnarray}
Then (\ref{isop1}) holds. Furthermore, if inequality 
(\ref{crucial2}) is strict, 
%or if $k > \frac{1}{3}$, 
then (\ref{M=BR}) holds only if $M = B_R $ for some $R > 0$. 
\end{theorem}

{\sl Proof:} We proceed similarly as in the proof of Theorem \ref{th1ter}. 
Below we mainly point out the differences, 
and we leave it to the reader to fill in the details.
\\ 
Note that our assumptions imply 
\begin{eqnarray}
\label{l/2}
 & & \frac{l}{2} <k \quad \mbox{and }
 \\
\label{positive} 
 & & 3k-2l-1 \geq 0 .
\end{eqnarray}
If $u\in C_0 ^{\infty } (\mathbb{R}^2 )\setminus \{ 0\}$, 
$u\geq 0$, we define $v$ by
$$
v(y):= u(x), \ \mbox{ where } \ y:= x|x| ^{\frac{3k-1}{4} } , \ \mbox{ and }
 \ z:= |y| .
$$
Then we show, using an interpolation argument 
as in the proof of Lemma \ref{4.3}, that for every $A\in 
\left[
0, \frac{16}{9(k+1)^2 }
\right] 
$, 
\begin{eqnarray}
\label{ineq1N=2}
 & & {\mathcal Q}_{k,l,2} (u) 
 \\
 \nonumber
 & \geq &
\left( 
\frac{3(k+1)}{4} 
\right) 
^{\frac{k+1}{l+2} } 
  \cdot
  \frac{
\left( 
\dint_{\mathbb{R} ^2 } |y|^{\frac{1}{3} } 
|\nabla _y v|\, dy 
\right) ^A  
\cdot 
\left( \dint_{\mathbb{R}^2 } |y|^{\frac{1}{3} } |v_z | \, dy 
\right) ^{1-A}
}{ 
\left( 
\dint_{\mathbb{R}^2 } |y|^{\frac{2(2l+1-3k)}{3(k+1)} } v^{\frac{l+2}{k+1}} \, dy 
\right)
^{\frac{k+1}{l+2} }
}
.
\end{eqnarray}
Let $\widetilde{v}$ denote the starshaped rearrangement of $v$.
Analogously as in the proof of Lemma \ref{4.2}, 
the properties of the rearrangement, (\ref{positive}) and Lemma \ref{Landes} lead to
\begin{eqnarray}
\label{starsh4N=2}
 & & 
 \int_{\mathbb{R}^2 } v^2 \, dy 
 = 
 \int_{\mathbb{R}^2 } 
 \widetilde{v} ^2  \, dy ,
 \\
 \label{starsh5N=2}
 & & 
 \int_{\mathbb{R}^2 } |y| ^{\frac{2(2l+1-3k)}{3(k+1)} } 
 v^{\frac{l+2}{k+1} } \, dy 
 \leq 
 \int_{\mathbb{R}^2 } |y| ^{\frac{2(2l+1-3k)}{3(k+1)} } 
 \widetilde{v} ^{\frac{l+2}{k+1} }\, dy,
 \\
 \label{vL1N=2}
 & & \frac{y \cdot \nabla \widetilde{v} 
 }{
 |y|^{\frac{2}{3} } }
 \equiv 
 |y|^{\frac{1}{3} } \frac{\partial \widetilde{v}
 }{
 \partial z} \in L^1 (\mathbb{R}^2 ), \ \mbox{ and}
 \\ 
 \label{starsh6N=2} 
 & & \int_{\mathbb{R}^2 } |y|^{\frac{1}{3} } 
 \left| 
 \frac{\partial v}{\partial z} 
 \right| \, dy 
 \geq 
 \int_{\mathbb{R}^2 } |y|^{\frac{1}{3} } 
 \left| 
 \frac{\partial \widetilde{v} 
 }{
 \partial z} \right| \, dy .
\end{eqnarray}
Further, we have by Hardy's inequality
\begin{equation}
\label{hardy2}
\int_{\mathbb{R}^2 } |y|^{\frac{1}{3} } |\widetilde{v}_z |\, dy 
\geq
\frac{4}{3} \int_{\mathbb{R}^2 } |y|^{-\frac{2}{3} } \widetilde{v} \, dy.
\end{equation}  
Finally, we show, similarly as in the proof of Lemma \ref{4.1}, 
that for every bounded measurable and starshaped set in $\mathbb{R} ^2 $,
\begin{eqnarray}
\label{holder2}
  & & \int_{M} |y| ^{\frac{4l-6k+2}{3k+3} } \, dy 
   \leq  d_2 \left( \int_M dy \right) ^{\frac{2(l+1-k)}{k+1} } \cdot 
  \left( \int_M |y|^{-\frac{2}{3} } \, dy \right) 
  ^{\frac{3k-2l-1}{k+1} }, 
  \quad \mbox{ where}
  \\
\label{d2}
 & & d_2 = \left( \frac{3}{2} \right) ^{\frac{2(l+1-k)}{k+1} } 
 \cdot \frac{k+1}{l+2} .
\end{eqnarray}  
By (\ref{crucial2}), we may choose 
$$
A= \frac{3(l+1-k)}{l+2} 
$$
in (\ref{ineq1N=2}). 
Combining this with (\ref{starsh5N=2}), (\ref{starsh6N=2}) and (\ref{hardy2}), 
we obtain
\begin{eqnarray}
\label{ineq1N=2new}
 & & {\mathcal Q}_{k,l,2} (u) 
 \\
\nonumber
 & \geq &
\left( \frac{3(k+1)}{4} 
\right) ^{\frac{k+1}{l+2} } 
  \cdot 
  \left( \frac{4}{3} \right) ^{\frac{3k-2l-1}{l+2} }
  \cdot 
  \frac{
\left( \dint_{\mathbb{R} ^2 } |y|^{\frac{1}{3} } |\nabla _y v|\, dy \right) 
^{\frac{3(l+1-k)}{l+2} }    
\cdot 
\left( \dint_{\mathbb{R}^2 } |y|^{- \frac{2}{3} } \widetilde{v} 
\, dy \right) ^{\frac{3k-2l-1}{l+2} }
}{ \left(
\dint_{\mathbb{R}^2 } |y|^{\frac{2(2l+1-3k)}{3(k+1)} } 
\widetilde{v}^{\frac{l+2}{k+1} } \, dy 
\right)  ^{\frac{k+1}{l+2} } }
.
\end{eqnarray}

Now let $M$ be a bounded measurable set,
and set
$M' := \{ y= x |x|^{\frac{3k-1}{4}} :\, x \in M\} $. 
Then, proceeding as the proof of Theorem
\ref{th1ter} and using the isoperimetric inequality 
\begin{equation}
\label{isoN=2}
{\mathcal R}_{\frac{1}{3},0,2} (M) \geq C_{\frac{1}{3},0,2} ^{rad} ,
\end{equation}
which follows from Corollary \ref{N=2_l<0}, 
we obtain from (\ref{ineq1N=2new}),
\begin{eqnarray}
\label{ineq6N=2}
 & & 
 {\mathcal R}_{k,l,2} (M )
 \\ 
 \nonumber
 & \geq &  
\left( \frac{3(k+1)}{4} 
\right) ^{\frac{k+1}{l+2} } 
  \cdot 
  \left( \frac{4}{3} \right) ^{\frac{3k-2l-1}{l+2} }
  \cdot 
  \frac{
\left( P_{\mu _{\frac{1}{3} }} (M')  \right) 
^{\frac{3(l+1-k)}{l+2} }    
\cdot 
\left( \dint_{\widetilde{M'} } |y|^{- \frac{2}{3} }  
\, dy \right) ^{\frac{3k-2l-1}{l+2} }
}{ \left(
\dint_{\widetilde{M'} } |y|^{\frac{2(2l+1-3k)}{3(k+1)} } 
 \, dy 
\right)  ^{\frac{k+1}{l+2} } }
\\
\nonumber
 & \geq &  
\left( \frac{3(k+1)}{4} 
\right) ^{\frac{k+1}{l+2} } 
  \cdot 
  \left( \frac{4}{3} \right) ^{\frac{3k-2l-1}{l+2} }
  \cdot 
  \frac{
\left( C_{\frac{1}{3},0,2} ^{rad} \mu _0  (M')  \right) 
^{\frac{3(l+1-k)}{l+2} }    
\cdot 
\left( \dint_{\widetilde{M'} } |y|^{- \frac{2}{3} }  
\, dy \right) ^{\frac{3k-2l-1}{l+2} }
}{ \left(
\dint_{\widetilde{M'} } |y|^{\frac{2(2l+1-3k)}{3(k+1)} } 
 \, dy 
\right)  ^{\frac{k+1}{l+2} } }
.
\end{eqnarray}
In view of (\ref{holder2}) and since $\mu_0 (M') =\mu_0 (\widetilde{M'})$, 
we finally obtain
\begin{equation}
\label{ineq7N=2} 
{\mathcal R}_{k,l,2} 
 \geq 
 \left( C_{\frac{1}{3} , 0,2} ^{rad} \right) ^{\frac{3(l+1-k)}{l+2} } 
 \cdot \left( \frac{3(k+1)}{4 d_2 } \right) ^{\frac{k+1}{l+2} } 
 \cdot \left( \frac{4}{3} \right) ^{\frac{3k-2l-1}{l+2} } = 
 C_{k,l,2} ^{rad} \ ,
\end{equation}
and (\ref{isop1}) follows by (\ref{CklNsmooth}). 
\\
Now assume that (\ref{M=BR}) holds. 
If inequality (\ref{crucial2}) is strict, 
then we must have $M=B_R $ for some $R>0$, by Lemma \ref{rangekl1}.
The Theorem is proved.
$\hfill \Box $
\begin{remark} 
\rm
The assumptions of Corollary \ref{N=2_l<0} and 
Theorem \ref{theoN=2k>0} cover the range $l=0\leq k\leq 1 $ for $N=2$.
\end{remark}
\noindent
\hspace*{1cm} 
Let us summarize the results of Section 5. 
\\
Let $k>1-N$ and $N\geq 2 $. We define a number $l_* = l_* (k,N)$ by  
\begin{equation}
\label{l*optimal}
l_* := \sup \{ l:\, l>-N,\, C_{k,l,N} = C_{k,l,N} ^{rad} \} .
\end{equation}
By Lemma \ref{rangekl1} we have that $C_{k,l,N} = C_{k,l,N} ^{rad} $ whenever 
$l\in (-N, l_* ]$. 
Further, Theorems \ref{R5}, \ref{th1bis} and Corollary \ref{sufficiency} tell us that
\begin{equation}
\label{l*k<0}
l_* = k\frac{N}{N-1} \quad \mbox{if }\ k\leq 0.
\end{equation}
Next, let $k> 0$ and define numbers $l^* = l^* (k,N) $ and $l_1 = l_1 (k,N) $ by
\begin{eqnarray}
\label{def_l^*}
l^* & := & k-1 + \frac{N-1}{k+N-1} ,
\\ 
\label{l_1N3}
l_1 (k,N ) & := & \frac{(k+N-1)^3 }{(k+N-1)^2 - \frac{(N-1)^2 }{N} } -N \quad \mbox{if }\ N\geq 3 ,
\\
\label{l_1N2}  
l_1 (k,2) & := &  \left\{ 
\begin{array}{ll}
0 & \mbox{ if $\ 0\leq k\leq \frac{1}{3}$}
\\
\frac{(k+1)^3 }{(k+1)^2 -\frac{16}{27}} -2 & \mbox{ if $\ k\geq \frac{1}{3} $}
\end{array}
\right.
.
\end{eqnarray} 
By Teorem \ref{R4} it follows that 
\begin{equation}
\label{l_*<l^*}
l_* \leq l^* .
\end{equation}
Further, Theorem \ref{th1ter}, Lemma \ref{csato} and 
Theorem \ref{theoN=2k>0} imply that 
\begin{equation}
\label{l_1<l^*}
l_1 \leq l^* .
\end{equation}
Note also that
\begin{equation}
\label{l_1<l_*}
l_1 \leq l_* .
\end{equation}
{\bf Conjecture 5.1 } {\sl There holds 
\begin{equation}
\label{conj1}
l_* (k,N)= l^* (k,N) \  \mbox{ when }\  k\geq 0.
\end{equation} 
}
Let us remark that equality \eqref{conj1} has already been conjectured in the case $N=2$ in \cite{DHHT}, Conjecture 4.22.

The approach used in the proof of  Theorem \ref{th1bis} also allows to obtain a lower 
bound for the isoperimetric constant $C_{k,l,N } $  for all positive values of $k$. Such a bound is useful when relation (\ref{isop1}) does not hold. In view of Theorem \ref{R4} this is the case when $l>l^* (k,N)$, or equivalently, if 
\begin{equation}
\label{violate}
l+1>k + \frac{N-1}{k+N-1} .
\end{equation}

\begin{proposition} 
\label{kpos}
Let $N\geq 2$, and assume $k\leq l+1$, $k> 0 $ and $l(N-1)/N \leq k$. Then 
\begin{equation} 
\label{mainpos} 
C_{k,l,N} \geq \left( \frac{N-1}{k+N-1} \right) ^{\frac{l+1-k}{l+N} } C_{0,l',N} ^{rad} ,
\end{equation}
where $l' :=  \frac{l(N-1)-kN}{k+N-1} $.
\end{proposition}

\begin{remark}\rm 
Similar estimates for the best constant $C_{k,l,N}$ have been obtained  in \cite{ChiHo}, Proposition 1.1, part 2, but with a different approach.
\end{remark}

{\sl Proof of Proposition 5.1 :\/} We proceed as in the proof of Theorem \ref{th1bis} until inequality \eqref{cambio2}. Then, since $\frac{N-1}{k+N-1}\leq 1$, we may replace  (\ref{cambio2}) by the 
inequality
\begin{equation}
\label{1stineq}
\int_{\mathbb{R} ^N}|x|^k |\nabla_x u| \, dx
\geq 
\frac {N-1}{k+N-1}
\int_{\mathbb{R} ^N} |\nabla_y v| \, dy .
\end{equation}
Continuing as before, we obtain
\begin{equation}
\label{2ndineq}
{\mathcal Q}_{k,l,N} (u) \geq \left( \frac{N-1}{k+N-1}\right) ^{\frac{l+1-k}{l+N}} {\mathcal Q}_{0,l', N} (v).
\end{equation}
Finally, observing that 
$l' = \frac{l(N-1)-kN}{k+N-1} \in [-1 ,0] $, (\ref{2ndineq}) yields (\ref{mainpos}). 
$\hfill \Box $

\noindent
\section{The case $N=1$}

The next result gives a complete solution to the  isoperimetric problem in the one-dimensional case.

\begin{theorem}  
\label{1d_Iso} 
Let $N=1$, $k>0 $ and $l>-1$.  
\\
{\bf (i)} If $k\geq l+ 1 $, then (\ref{isop1}) holds. 
Moreover, if (\ref{M=BR}) holds and if $k>l+1$, then $M= (-R,R)$ for some $R>0$.    
\\
{\bf (ii)} If $k<l+1$, then 
\begin{eqnarray}
\label{N=1_2}
 & & {\mathcal R}_{k,l,1} (M ) \geq  {\mathcal R}_{k,l,1} ((0,1))  
= (l+1) ^{k/(l+1)}  \\
\nonumber
 & &  \mbox{for all measurable sets $M$ with $0<\mu_l (M) <+\infty $.} 
\end{eqnarray}
\end{theorem}

{\sl Proof: }
\\
{\bf (i) }  The result follows from Theorem \ref{R5}.
\\
{\bf (ii) } It is sufficient to prove the assertion for smooth sets, that is, for unions of finitely many bounded open intervals. 
For any smooth set $\Omega $ we set
$$
U:= \{ y= |x|^l x :\, x\in \Omega \} .
$$
Then an elementary calculation shows that
\begin{equation}
\label{N=1_3}
{\mathcal R}_{k,l,1} (\Omega ) = (l+1) ^{k'} {\mathcal R}_{ k', 0,1} (U),
\end{equation}
where $k' = \frac{k}{l+1} \in (0,1)$.
It remains to show that
\begin{equation}
\label{N=1_4}
{\mathcal R}_{ k', 0,1} (U)\geq {\mathcal R}_{ k', 0,1} ((0,1)) \quad 
\mbox{for all smooth sets $U\subset \mathbb{R} $. }
\end{equation} 
Let $y_1 := \inf U$ and $y_2 :=\sup U$. 
Then $|y_1 | ^{k'} + |y_2 | ^{k'} = P_{\mu_{k'}} ((y_1 ,y_2 )) 
\leq P_{\mu _{k'} } (U) $ and $\int _U \, dy   \leq 
\int_{y_1 } ^{y_2 } \, dy $. In other words, we have
$$
{\mathcal R}_{k',0,1} (U)\geq {\mathcal R}_{k',0,1} ((y_1 , y_2 )).
$$ 
It is therefore sufficient to consider open intervals $U$.
Thus, let  
$U= (y_1 , y_2 )$, ($y_1 <y_2 $).
Setting $c:= y_2 -y_1 $, we define
$$
U(t) := (-c/2 +t,  c/2 +t), \quad (t\in \mathbb{R} ).
$$
Then we have $\int _{U(t)} dy=c $ and 
$$
\int_{\partial U(t)} |y| ^{k'} \, \mathscr{H}_0 (dy) = 
|-c/2 +t| ^{k/(l+1)} + |c/2 +t| ^{k/(l+1)} =:f(t), \quad (t\in \mathbb{R} ).  
$$
Note that $f$ is an even function.
Let $t\in [-c/2, c/2] $. Then $f(t) = (c/2 -t) ^{k/(l+1)} + 
(t+ c/2 ) ^{k/(l+1)} $, which is a concave function. Hence 
$$
\inf \{ f(t):\, t\in [-c/2 , c/2 ] \} = f(-c/2 )= f(c/2).
$$
Since also $f'(t)>0$ for $t>c/2$, this implies that  
$$
\inf \{ f(t):\, t\in \mathbb{R} \} = f(-c/2 )= f(c/2),
$$  
that is,
$$
{\mathcal R}_{k',0,1} ((y_1 , y_2 )) \geq {\mathcal R}_{k',0,1} ((0,c)),
$$
and the assertion follows. 
$\hfill \Box $
\\[0.1cm]
\medskip

\section{The case $l+N<0$}
In this section we treat our functionals 
${\mathcal R}_{k,l,N} $ and ${\mathcal Q}_{k,l,N} $ 
for a different range of the parameters $k$ and $l$. 
Instead of (\ref{ass1}) we assume
\begin{equation}
\label{ass2}
k+N-1 <0 \quad \mbox{and }\ l+N<0 .
\end{equation}

We state our result only for smooth sets. 
Extensions to measurable sets and a discussion 
of the equality case in the isoperimetric inequalities 
follows the lines of the proofs in section 4, and they are left to the reader.  
	
\begin{theorem}
\label{maintheorem2}
Let $N\in \mathbb{N} $, 
$k,l\in \mathbb{R} $ and $l+N <0$. Further, assume that one of the following conditions holds:
\\
{\bf (j)} $N\geq 1 $ and $l+1\geq k $;
\\
{\bf (jj)} $N\geq 2$, $l+1\leq k$, $ k\leq l\frac{N-1}{N} $ and $k +2N-2 \geq 0 $; 
\\ 
{\bf (jjj)} $N\geq 3$, $ l+1\leq k\leq 2-2N $ and 
$$
\frac{1}{l+N } \leq \frac{1}{k+N-1} - \frac{ (N-1)^2 }{ N(k+N-1)^3 }.
$$ 
{\bf (jv)} $N= 2 $, $l+1\leq k$, $\frac{l}{2} \leq k$ and either 
\begin{eqnarray*}
 & & -2\leq k\leq -\frac{7}{3}  \ \mbox{ or }
 \\
 & & -\frac{7}{3} \leq k \ \mbox{ and } \
 \frac{1}{l+2} \leq \frac{1}{k+1} - \frac{16}{27 (k+1)^3 } ;
\end{eqnarray*} 
Then  
\begin{equation}
\label{mainineq2}
\dint_{\partial \Omega } |x|^k \, \mathscr{H}_{N-1} (dx)
\geq 
\overline{C}_{k,l,N} ^{rad}  
\left( 
\int_{\Omega } |x|^l \, dx 
\right) 
^{(k+N-1)/(l+N) } , 
\end{equation}
for every open set $\Omega \subset \mathbb{R}^N $ with smooth boundary that does 
not contain a neighborhood of the origin,
where 
\begin{equation}
\label{Cklinverse}
\overline{C}_{k,l,N} ^{rad} := (N\omega _N ) ^{(l-k+1)/(l+N)} \cdot  |l+N | ^{(k+N-1)/(l+N)} .
\end{equation}
Equality in (\ref{mainineq2}) holds  for all sets $\Omega = \mathbb{R} ^N \setminus \overline{B_R }$, ($R>0$).
\end{theorem}

\begin{remark}
\rm
Theorem 7.1  has been known in some particular situations:  
\\
{\bf 1.} $N=2 $, $k=l<-2 $, see \cite{CJQW}, Proposition 4.3;
\\
{\bf 2.} $N\in \mathbb{N} $, $k=l<-N $, see \cite{DHHT}, Proposition 7.5;
\\
{\bf 3.} case {\bf (jj)}, see \cite{ChiHo}, Theorem 1.3, part 3. 
\end{remark}
{\sl Proof of Theorem 7.1 : } 
 Let $u\in C^{\infty } _0(\mathbb{R}^N)\setminus \{ 0\} $, with 
$u \not\equiv  0$  in $\mathbb{R} ^N$.
We set
 $$
 y:=x|x|^{-2}\, , \quad v(y):=u(x)\,
 .
 $$
 Observe that $v$ vanishes in a neighborhood of the origin. 
Then a short computation shows that
\begin{eqnarray}
 \label{uext}
\int_{\mathbb{R} ^N}|x|^l
|u|^{(l+N)/(k+N-1)} \, dx
 & = &
\int_{\mathbb{R} ^N}
|y|^{-l-2N}
 |v|^{(l+N)/(k+N-1)}\, dy
 \quad \mbox{and }
\\
\int_{\mathbb{R} ^N}|x|^k |\nabla_x u| \, dx
 & = & 
\int_{\mathbb{R} ^N}|y|^{-k-2N+2} |\nabla_y v| \, dy.
\end{eqnarray}
This implies that
\begin{eqnarray}
\label{QQ}
 & & Q_{k,l,N}(u)  =  Q_{\widetilde{k}, \widetilde{l},N }(v),
\\
\nonumber 
 & & \mbox{where  $\widetilde{k} := -k-2N+2$ and $\widetilde{l} := -l-2N$.} 
\end{eqnarray}
(\ref{QQ}) also means that for every open set $\Omega $ with smooth boundary that does not contain a neighborhood of the origin,
\begin{equation}
\label{RR}
 R_{k,l,N}(\Omega )  =  R_{\widetilde{k}, \widetilde{l},N }(\widetilde{\Omega }), 
\quad 
\mbox{where  $\widetilde{\Omega } := \{ y= \frac{x}{|x|^2 } : \, x\in \Omega \} .$}
\end{equation}
Now the conclusion follows from Theorem \ref{maintheorem}. 
$\hfill \Box $

\section{Applications}
In this section we provide some applications of our results.
\subsection{Polya-Szeg\"o principle}
First we obtain a Polya-Sz\"ego principle related to our 
isoperimetric inequality (\ref{isop1}).
Assume that the numbers $l$ and $k$ satisfy one of the conditions {\bf (i)}-{\bf (iv)} of 
Theorem \ref{maintheorem}. Then (\ref{mainineq}) implies
\begin{equation} 
\label{Isop_k_l}
\int_{\partial \Omega }
|x|^k {\mathscr H}_{N-1}(dx)
\geq 
\int_{\partial \Omega ^{ \star }}
|x|^k 
{\mathscr H}
_{N-1}(dx)
\end{equation}
for every smooth set $\Omega $, where $\Omega ^{\star}$ is the $\mu_l $-symmetrization of $\Omega $.
We will use (\ref{Isop_k_l}) to prove the following

\begin{theorem} 
\label{ps}
(Polya-Szeg\"o principle) 
Let the numbers $k,l$ and $N$ satisfy one of the conditions {\bf (i)}-{\bf (iv)} 
of Theorem \ref{maintheorem}.
Further, let $p\in [1, +\infty)$ and $m:= pk+(1-p)  l $. Then there holds 
\begin{equation}
\int_{\mathbb{R}^{N}}\left\vert \nabla u\right\vert ^{p}\left\vert
x\right\vert ^{pk+(1-p)l}dx\geq \int_{\mathbb{R}^{N}}\left\vert \nabla
u^{ \star }\right\vert ^{p}\left\vert x\right\vert ^{pk+(1-p)l}dx \quad 
\forall u\in W_0 ^{1,p} (\mathbb{R} ^N , d\mu _m ),
\label{PS_k_l}
\end{equation}
where $u^{\star } $ denotes the $\mu _l $-symmetrization of $u$.
\end{theorem}

\noindent {\sl Proof:} 
It is sufficient to consider the case that $u$ is non-negative. Further,
by an approximation argument we may assume that
$u \in C^{\infty}_{0}(\mathbb{R}^{N} )\setminus \{ 0 \} $.
Let 
\begin{eqnarray*}
I & := &  
\int_{\mathbb{R}^{N}} | \nabla u| ^{p}
|x| ^{pk+(1-p)l} \, dx \quad \mbox{and}\\
I ^{\star }  & := &  
\int_{\mathbb{R}^{N}} | \nabla u^{\star} | ^{p}
|x| ^{pk+(1-p)l} \, dx .
\end{eqnarray*}
The coarea formula yields 
\begin{eqnarray}
\label{1coarea}
I 
 & = & 
 \int_{0}^{\infty }\int_{u=t} |\nabla
u| ^{p-1} |x| ^{pk+(1-p)l} \, {\mathscr H}_{N-1}(dx)\, dt \quad \mbox{and}
\\
\label{coarea2}
I^{\star}
 & = & 
 \int_{0}^{\infty }\int_{u^{\star } =t} |\nabla
u^{\star} | ^{p-1} |x| ^{pk+(1-p)l} \, {\mathscr H}_{N-1}(dx)\, dt 
.
\end{eqnarray}
Further, H\"older's inequality gives
\begin{equation}
\label{1holder}
\int_{ u=t} |x|^k \, \mathscr{H} _{N-1} (dx)
\leq 
\left( \int_{ u=t} |x|^{kp +l(1-p)} |\nabla u| ^{p-1} \, \mathscr{H} _{N-1} (dx) \right) ^{\frac{1}{p} } 
\cdot
\left( \int_{ u=t} \frac{|x|^l }{|\nabla u| } \, \mathscr{H} _{N-1} (dx) \right) ^{\frac{p-1}{p} } ,
\end{equation} 
for a.e. $t\in [0, +\infty )$. 
Hence (\ref{1coarea}) together with (\ref{1holder}) tells us that
\begin{equation}
\label{coarea3}
I 
\geq 
\int_{0}^{\infty }
\left( \int_{u=t} |x| ^{k} \, {\mathscr H}
_{N-1}(dx)
\right) ^{p} \cdot 
\left( 
\int_{u=t}\frac{ |x| ^{l}}{
| \nabla u| } \, {\mathscr H}_{N-1}(dx)
\right) ^{1-p} \, dt.
\end{equation}
Since $u^{\star} $ is a radial function, we obtain in an analogous manner,
\begin{equation}
\label{coarea4}
I^{\star} 
=
\int_{0}^{\infty }
\left( \int_{u^{\star} =t} |x| ^{k} \, {\mathscr H}
_{N-1}(dx)
\right) ^{p} \cdot 
\left( 
\int_{u^{\star} =t}\frac{ |x| ^{l}}{
| \nabla u^{\star} | } \, {\mathscr H}_{N-1}(dx)
\right) ^{1-p} \, dt.
\end{equation}
Observing that
\begin{equation}
\label{meas_u>t}
\int_{u>t} |x|^{l} \, dx
=
\int_{u^{\star }>t}
|x|^{l} \, dx \quad  \forall t\in [0, +\infty ), 
\end{equation}
Fleming-Rishel's formula yields
\begin{equation}
\label{flemingrishel} 
\int_{u=t } \frac{|x|^l }{|\nabla u|} \, \mathscr{H}_{N-1} (dx) 
=
\int_{u^{\star} =t } \frac{|x|^l }{|\nabla u^{\star} |} \, \mathscr{H}_{N-1} (dx)
\end{equation}
for a.e. $t\in [0, +\infty )$. 
Hence
(\ref{flemingrishel}) and (\ref{Isop_k_l}) give
\begin{eqnarray*}
 & &
\int_{0}^{\infty }
\left( \int_{u=t} |x|^k \, \mathscr{H}
_{N-1}(dx) \right) ^{p}
\cdot 
\left( \int_{u=t}\frac{| x| ^{l}}{
| \nabla u| } \, \mathscr{H}_{N-1}(dx) \right) ^{1-p}
\, dt 
\\
 & \geq &
\int_{0}^{\infty }\left( \int_{u^{\star} =t} |x| ^{k}
\, \mathscr{H}_{N-1}(dx) \right) ^{p} \cdot \left( \int_{u^{\star}=t}
\frac{|x|^{l}}{|\nabla u^{\star} | } \, \mathscr{H}
_{N-1}(dx)\right) ^{1-p} \, dt.
\end{eqnarray*}
Now (\ref{PS_k_l})  follows from this, (\ref{coarea3}) and (\ref{coarea4}).
$\hfill \Box$
\\[0.1cm]
An important particular case of Theorem \ref{ps} is 
\begin{corollary}
\label{specialcasePS}
Let $p\in [1, +\infty )$, $a\geq 0 $, $u\in W_0 ^{1,p} 
(\mathbb{R}^N , d\mu _{ap }) $, and let $u^{\star } $ be the Schwarz 
symmetrization (=$\mu_0 $-symmetrization) of $u$.
Then
\begin{equation}
\label{PSspecial}
\int_{\mathbb{R}^N } \left| \nabla u\right|^p  |x|^{ap} \, dx 
\geq 
\int_{\mathbb{R}^N } \left| \nabla u^{\star} \right|^p |x|^{ap} \, dx.
\end{equation}
\end{corollary}

{\sl Proof: } We choose $k:= a $ and $l:= 0$. If $a\in [0,1]$ then $k,l$ 
satisfy either one of the conditions {\bf (iii)} or {\bf (iv)}, and if $a\geq 1 $, then $k,l$ satisfy condition {\bf (i)} of Theorem \ref{maintheorem}. Hence (\ref{PSspecial}) follows from Theorem \ref{ps}.
$\hfill \Box $
   
\subsection{Caffarelli-Kohn-Nirenberg inequalities}
Next we will use Theorem \ref{ps} to obtain best constants in some 
Caffarelli-Kohn-Nirenberg inequalities. 

Let $p,q, a, b$ be real numbers  such that
\begin{eqnarray}
\label{CKNassump1}
 & & 1\leq p \leq q \left\{ 
 \begin{array}{ll} 
 \leq \frac{Np}{N-p} &  \mbox{ if } \ p< N
 \\
 < +\infty & \mbox{ if } \ p\geq N 
 \end{array}
 \right.
 , 
 \\
\nonumber
 & & a> 1-\frac{N}{p}, \quad \mbox{and }
 \\
\nonumber
 & & b= b(a,p,q,N) = N \left( \frac{1}{p} -\frac{1}{q} \right) + a-1 .
\end{eqnarray}
We define
\begin{eqnarray}
\label{p*}
p^* & := & \left\{ 
\begin{array}{ll} 
\frac{Np}{N-p} & \mbox{ if } p<N
\\
+\infty & \mbox{ if } p\geq N
\end{array}
\right.
,
\\
& &\nonumber
 \\
\label{fctalE}
E_{a,p,q,N} (v)
 & := & 
\frac{\dint_{\mathbb{R} ^N } |x|^{ap} |\nabla v|^p \, dx
 }{
\left( \dint_{\mathbb{R}^N } |x|^{bq} |v|^q \, dx \right) ^{p/q}  }, 
\quad v\in C_0 ^{\infty } (\mathbb{R}^N  )\setminus \{ 0\}  ,
\\
& &\nonumber
 \\
S_{a,p,q,N} & := & \inf \{ E_{a,p,q,N} (v): \, v\in C_0 ^{\infty } 
(\mathbb{R}^N  ) \setminus \{ 0\} \}, \quad \mbox{and}
\\
& &\nonumber
 \\
S_{a,p,q,N} ^{rad} & := & \inf \{ E_{a,p,q,N} (v): \, v\in C_0 ^{\infty } 
(\mathbb{R}^N  )\setminus \{ 0\} ,  \ v \mbox{ radial }\}.
\end{eqnarray}
It has been proved in \cite{CKN}, that 
\begin{equation}
\label{S>0}
S_{a,p,q,N}>0 .
\end{equation}
Further, it is known that the functional $E_{a,p,q,N} $ is well-defined for functions in $W_{0} ^{1,p} (\mathbb{R}^N , d\mu _{ap} )$ and that  
$C_0 ^{\infty } (\mathbb{R} ^N )$ is dense in 
$W_0 ^{1,p} (\mathbb{R} ^N , d\mu _{ap} )$. 
Moreover,  
$S_{a,p,q,N}$ is attained for some $u\in W_0 ^{1,p} (\mathbb{R} ^N , d\mu _{ap } )$ if
$1<p<q $.
\\
We are interested 
in the range of values $a$ (depending on $p,q$ and $N$) for which 
\begin{equation}
\label{S=S_rad} 
S_{a,p,q,N} = S_{a,p,q,N} ^{rad}
\end{equation} 
holds. 
This problem has been investigated by several authors. For recent advances concerning the symmetry of optimizers in the CKN
inequalities, see for example \cite{LL},  \cite{DolEstLoss} and references therein.
\\
First observe that the case $1<p=q$ (which is equivalent to $a-b=1$) 
corresponds to the Hardy-Sobolev inequality, with the known best constant
\begin{equation}
\label{hardyconstant}
S_{a,p,p,N} = S_{a,p,p,N} ^{rad} =
\left( \frac{N}{p} -1+a\right) ^p ,
\end{equation}
see \cite{Hardy}. Note that the Hardy constant $S_{a,p,p,N} $ is not achieved for any function
$u\in W_0 ^{1,p} (\mathbb{R} ^N , d\mu _{ap} )$. 
\\ 
Next, let $1<p <N$ and $q=p^* $. If $a\leq 0$, then one has 
\begin{equation}
\label{s=sradcrit}
S_{a,p,p^* ,N} = S_{a,p,p^* ,N} ^{rad} ,
\end{equation}
see \cite{HK}, 
Theorem 2.4, condition (3).
\\
\hspace*{1cm} From now on let us assume that
\begin{equation}
\label{q<p*}
N\geq 2 \ \mbox{ and }\ 
1<p<q<p^* .
\end{equation} 
In this case, the constants $S_{a,p,q,N} ^{rad} $, including the corresponding (radial) minimizers, 
have been given in \cite{Musina}, Theorem 1.4.  The problem of symmetry breaking was analyzed by many authors, see \cite{CaldMus} and the references cited therein. 
\\
It is known that there is a finite number 
$$
a_* = a_* (p,q,N) 
$$
with $a_* \geq 1-\frac{N}{p} $, such that 
\begin{eqnarray}
\label{a<=a_*} 
S_{a,p,q,N} = S_{a,p,q,N} ^{rad} & & \mbox{for }\ a\in ( 1-\frac{N}{p} ,a_*] \ \mbox{ and }
\\
\label{a>a_*} 
S_{a,p,q,N} < S_{a,p,q,N} ^{rad} & & \mbox{for }\ a> a_* ,
\end{eqnarray}
see \cite{CaldMus}, Theorem 1.1 and Remark 3.1.
Moreover, if $a^* = a^* (p,q,N)$ denotes the unique number in $( 1-\frac{N}{p} , +\infty )$ such that
\begin{equation} 
\label{def_a^*}
 \left( \frac{N}{p} -1 +a^*  \right)
^2 = (N-1) \left( \frac{1}{q-p} -\frac{1}{q+p'} \right) ,
\end{equation}  
where $p' = \frac{p}{p-1} $, then 
\begin{equation}
\label{a_*<a^*}
a_* \leq a^* ,
\end{equation} 
see \cite{CaldMus}, Theorem 1.1, and if 
 $p<N$, then $a_* \geq 0$, see \cite{CaldMus}, Theorem 1.3. 
\\  
Finally, it has been conjectured that condition (\ref{a_*<a^*}) cannot be improved, see  \cite{CaldMus}, p. 423, that is: 
\\[0.1cm]
{\bf Conjecture 8.1 } {\sl There holds 
\begin{equation}
\label{conj2} a^* =a_* .
\end{equation}
} 
\begin{remark}
\rm 
The case $p=2$ in the CKN inequalities has received a lot of interest since the  seminal article \cite{CW}. In particular, Conjecture 8.1 for $p=2$ has been proved in the recent paper \cite{DolEstLoss}, Theorem 1.1, using generalized entropy functionals for diffusion equations.  However, this tool seems not useful for general $p$.
\\[0.1cm]
\hspace{0.3cm}
A bound  for $a_* $ from below is given in \cite{HK}, 
Proposition 4.6: Let
\begin{equation}
\label{a_1}
a_1 = a_1 (p,q,N) :=  \frac{N-1}{1+ \frac{q}{p'} } - \frac{N}{p} +1 ,
\end{equation}
and note that $a_1 >1- \frac{N}{p} $.
Then
\begin{equation}
\label{a*boundbelow}
a_* \geq a_1.
\end{equation} 
\end{remark}
Our aim is to improve on the bound $a_1 $. 
First observe that an application of the Theorems \ref{maintheorem} and \ref{ps} yield the following result. 
\begin{lemma}
\label{CKN}
Assume that $N, p,q,a$ and $b$ satisfy the conditions (\ref{CKNassump1}) and (\ref{q<p*}).  
Further, assume that there exist real numbers $k$ and $l$ which satisfy one of the conditions 
{\bf (i)}-{\bf (iv)} of Theorem \ref{maintheorem}, and such that
\begin{eqnarray}
\label{akl}
 & & ap = kp + l(1-p) \ \mbox{ and }
 \\
\label{bq<l}
 & & bq \leq l.
\end{eqnarray} 
Then (\ref{S=S_rad}) holds.
\end{lemma}

{\sl Proof:} Let $u\in W_0 ^{1,p} (\mathbb{R} ^N , d\mu_{ap} )\setminus \{ 0\} $,  
and let $u^{\star} $ be the $\mu_l $-symmetrization of $u$. 
Then we have by Theorem \ref{ps} and (\ref{akl}), 
\begin{equation}
\label{ps1}
\int_{\mathbb{R} ^N } |x|^{ap} |\nabla u| ^p \, dx \geq  
\int_{\mathbb{R} ^N } |x|^{ap} |\nabla u^{\star}| ^p \, dx.
\end{equation}
Further, it follows from (\ref{hardylitt1}) and (\ref{bq<l})
that 
\begin{equation}
\label{bqint}
\int_{\mathbb{R} ^N } |x|^{bq} | u| ^q \, dx \leq  
\int_{\mathbb{R} ^N } |x|^{bq} | u^{\star}| ^q \, dx.
\end{equation}
Finally, (\ref{ps1}) together with (\ref{bqint}) yield
\begin{equation}
\label{E>E*}
E_{a,p,q,N} (u) \geq E_{a,p,q,N} (u^{\star} ),
\end{equation}
and the assertion follows.
$\hfill \Box $
\\[0.1cm]
Next we define 
\begin{equation}
\label{a2}
a_2 = a_2 (p,q,N) :=  1+ N\left( \frac{1}{q} -\frac{1}{p} \right) ,
\end{equation}
and note that 
\begin{equation}
\label{a2a1}
\max \{ 0, a_1 \} <a_2 <1.
\end{equation}
\begin{proposition}
\label{improve1}
Assume that $N, p,q,a$ and $b$ satisfy the conditions 
(\ref{CKNassump1}) and (\ref{q<p*}), 
and let 
\begin{equation}
 a\leq a_2 .
\end{equation}
Then (\ref{S=S_rad}) holds.
\end{proposition}
{\sl Proof :} We may restrict to the case $a\geq 0$, and we choose 
 $k:=a$ and $l:=0$. 
Since $0< k< 1$, 
one of the conditions {\bf (iii)} or {\bf (iv}) of 
Theorem \ref{maintheorem} is satisfied. 
Further, we have 
\begin{eqnarray*}
bq-l & = & bq = \left( N\left( \frac{1}{p} -\frac{1}{q} \right) +a-1\right) q 
\\
 & \leq &  \left( N\left( \frac{1}{p} -\frac{1}{q} \right) +a_2 -1\right) q =0.
 \end{eqnarray*}
Now the assertion follows from Lemma \ref{CKN}.
$\hfill \Box $
\\[0.1cm]
Finally, a more sophisticated choice of the parameters $k$ and $l$ leads to a 
further improvement of the lower bound for $a_* $.
\\
First we assume $N\geq 3$. Let
us define $a_3 = a_3 (p,q,N)$ as the unique number in $\left( 1- \frac{N}{p} , +\infty \right) $, 
such that
\begin{equation}
\label{def_a_3}
\left( \frac{N}{p} -1 +a_3 \right) ^2 = 
\frac{(N-1)^2 }{ N\left( \frac{1}{p}- \frac{1}{q} \right) \cdot
\left( 1-\frac{q}{p} +q \right) ^2 }
.
\end{equation}
Note that
\begin{equation}
\label{a_2<a_3 }
a_2 < a_3 .
\end{equation}
\begin{theorem}
\label{2ndimprove}
Assume that $N, p,q,a$ and $b$ satisfy $1<p<q<p^*$, 
$N\geq 3$ and the conditions (\ref{CKNassump1}). 
Further, let
\begin{equation}
 a\leq a_3 .
\end{equation}
Then (\ref{S=S_rad}) holds.
\end{theorem}
{\sl Proof :}   
By elementary calculus one verifies that $a_3 $ appears as the maximum 
of all values $a\geq 0$ which have a representation
$a= k + l(\frac{1}{p} -1) $ with parameters $k$ and $l$ 
that satisfy the conditions {\bf (iii)} of Theorem \ref{maintheorem} 
and such that $bq\leq l $. Formally,   
\begin{eqnarray}
\label{optimize_a}
 a_3 & = & \max \left\{ a: \, a= k+l(\frac{1}{p} -1) , \, 0\leq k\leq l+1, \right. 
\\
 \nonumber 
 & & \qquad \left. \frac{1}{l+N} \geq \frac{1}{k+N-1} - 
\frac{(N-1)^2 }{N (k+N-1) ^3 }  ,
\, bq\leq l \right\} .
\end{eqnarray}
The assertion now follows from Lemma \ref{CKN}.
$\hfill \Box $  
\\[0.1cm]
The bound $a_2 $ can be improved in the case $N=2$, too, provided that
\begin{equation}
\label{pqspecial}
\frac{1}{q} > \frac{1}{p} -\frac{1}{3} .
\end{equation}
Define $a_4 = a_4 (p,q) $ as the unique number in $(1-\frac{2}{p} , + \infty )$ 
such that
\begin{equation}
\label{def_a_4}
\left( \frac{2}{p} -1+ a_4 \right) ^2 = 
\frac{16}{27 \left( \frac{1}{p} -\frac{1}{q} \right) \left( 1-\frac{q}{p} +q \right) ^2 } .
\end{equation}
Note that
\begin{equation}
a_2 (p,q,2) = 1+ 2\left( \frac{1}{q} -\frac{1}{p} \right)  < a_4 ,
\end{equation} 
in view of (\ref{CKNassump1}) and (\ref{pqspecial}).
\begin{theorem}
\label{3rdimprove}
Let $N=2$, and assume that the numbers $N,p,q,a$ and $b$ satisfy $1<p<q<p^* $ and the conditions  
(\ref{CKNassump1}), (\ref{pqspecial}). 
Further, let
\begin{equation}
 a\leq a_4 .
\end{equation}
Then (\ref{S=S_rad}) holds.
\end{theorem}
{\sl Proof :}   Using the conditions {\bf (iv)} of Theorem \ref{maintheorem}
one verifies that    
\begin{eqnarray}
\label{optimize_a_4}
 a_4 & = & \max \left\{ a: \, a= k+l(\frac{1}{p} -1) , \, \frac{1}{3}\leq k\leq l+1, \right. 
\\
 \nonumber 
 & & \qquad \left. \frac{1}{l+2} \geq \frac{1}{k+1} - \frac{ 16}{27 (k+1)^3} ,\, bq\leq l \right\} .
\end{eqnarray}
Note that the set on the right-hand side of (\ref{optimize_a_4}) 
is non-empty in view of (\ref{pqspecial}).
Now  the assertion again follows from Lemma \ref{CKN}.
$\hfill \Box $

\begin{remark}
\rm 
Let us point out an interesting relation between the two numbers $l^* $ and $a^* $ defined by  
(\ref{def_l^*}) respectively (\ref{def_a^*}).
\\
First observe that condition (\ref{def_a^*}) appears as a limit case of (\ref{def_l^*}) 
when sending $p\to 1 $ and putting
$$
a^* := k,\ q:= \frac{l^* +N}{k+N-1} , \ \mbox{ and } b:=\frac{l^* (k+N-1)}{l^* +N} .
$$
Further, assume that Conjecture 5.1 was true. Then, proceeding similarly as in the proof of Theorem \ref{2ndimprove}, one can show that also Conjecture 8.1 holds true: Indeed,  by elementary calculus one verifies that 
\begin{equation}
\label{optimize_a^*}
a^*  =  \max \left\{ a: \, a= k+l(\frac{1}{p} -1) , \ 0\leq k\leq l+1\leq k + \frac{N-1}{k+N-1} ,  
\ bq\leq l \right\} .
\end{equation}  
Then one obtains as before that 
$$
E_{a,p,q,N} (u) \geq E_{a,p,q,N} (u^{\star} ) 
\qquad \forall u\in W_0 ^{1,p} (\mathbb{R}^N , d\mu _{ap} )\setminus \{ 0\} ,
$$  
and (\ref{conj2}) follows.   
\end{remark}   

\subsection{Sobolev-type inequalities for Lorentz spaces} 
Corollary \ref{specialcasePS} can be used to obtain 
best constants for imbedding inequalities between the 
Sobolev space $W_0 ^{1,p} (\mathbb{R} ^N , d\mu _{ap } )$, 
with $a\geq 0 $, into Lorentz spaces. 
\\  
Let  $u: \mathbb{R}^N \to \mathbb{R} $ be a measurable function and 
$u^{ \star} $ its Schwarz symmetrization ($=\mu_0 $-symmetrization). 
Then the  decreasing rearrangement  of $u$ is given by  
\begin{equation}
\label{schwarz_decr}
u^* (\omega _N |x|^N ) = u^{\star} (x ), \quad (x\in \mathbb{R}^N ).
\end{equation}
For every $r\in (0,\infty )$ we define
\begin{eqnarray}
\label{seminorma}
\Vert u\Vert _{r,q}
 & = &
\left(
\int_0 ^{+\infty}
\left[ u^*  (s)\, s^{1/r}\right] ^q
\, \frac{ds}{s}
\right)
^{1/q}\,
 \quad \mbox{if $q\in (0,\infty ) $,  and }
\\ 
\label{marc}
 \Vert u \Vert _{r,\infty}
 & = &
\sup _{s>0}\, u^{ \star} (s)\, s^{1/r},
\quad \mbox{if $q=+\infty$.}
\end{eqnarray} 
The {\sl Lorentz space\/}  $L^{r,q}(\mathbb{R} ^N )$ is the collection of all  
measurable functions $u: \mathbb{R}^N \to \mathbb{R}$ such that 
$\Vert u \Vert _{r,q}$ is finite.
These spaces give in some sense a refinement of the usual Lebesgue spaces. 

\begin{theorem}
\label{lorentz}
Let $N,a,p,q$ and $b$ satisfy the conditions (\ref{CKNassump1}) and (\ref{q<p*}), with $a\in 
[0,  a_2 ]$, 
where $a_2 $ is given by (\ref{a2}). 
Then we have 
\begin{equation}
\left( \int_{\mathbb{R} ^N } |x|^{ap} |\nabla u| ^p \, dx \right) ^{1/p} 
\geq (\omega _N )^{-b/N} 
\left( S_{a,p,q,N} ^{rad} \right) ^{1/p} \Vert u\Vert _{r,q} \qquad \forall u\in W_0 ^{1,p} (\mathbb{R} ^N , d\mu _{ap} ) ,
\end{equation}
where
\begin{equation}
\label{defi_r}
r := \frac{Np}{N-p+ap}.
\end{equation}
\end{theorem}

{\sl Proof:} Let $u\in W_0 ^{1,p} (\mathbb{R} ^N , d\mu _{ap} )
\setminus \{ 0\} $, and let $u^{\star} $ denote its Schwarz symmetrization.  
Corollary  \ref{specialcasePS} tells us that  
\begin{equation}
\label{ps2}
\int_{\mathbb{R} ^N }  |\nabla u| ^p |x|^{ap}\, dx 
\geq 
\int_{\mathbb{R} ^N } |\nabla u^{\star} | ^p |x|^{ap}\, dx .
\end{equation} 
Further, 
%since $a\le a_2$, we get
%$$
%b:= N\left( \frac{1}{p} -\frac{1}{q} \right) -1+a \leq 0,
%$$
we have by 
Proposition \ref{improve1},  
\begin{equation}
\label{comp_u*}
\int_{\mathbb{R} ^N } |\nabla u^{\star} |^p |x|^{ap} \, dx \geq S_{a,p,q,N} ^{rad} 
\left( \int_{\mathbb{R}^N } |u^{\star} |^q |x|^{bq} \, dx \right) ^{p/q} .
\end{equation} 
Since also
$$
\Vert u\Vert _{r,q} =  \left( \omega _N \right) ^{-b/N}
\cdot
\left( \int_{\mathbb{R}^N } |u^{\star} |^q |x|^{bq} \, dx \right) ^{1/q} ,
$$ 
where $r$ is given by (\ref{defi_r}),
the assertion follows from (\ref{ps2}) and (\ref{comp_u*}). 
$\hfill \Box $

\begin{remark}
\rm
{\bf (a)} 
Theorem \ref{lorentz} is well-known in the special case $a=0$, see
\cite{Alvino}, where also other cases are considered, and \cite{E}  
and \cite{CRT}.
\\
{\bf (b)} Note that the number $r$ defined in (\ref{defi_r}) satisfies
\begin{equation}
\label{range_r}  
q\leq r<+\infty,
\end{equation}
by the assumptions of Theorem \ref{lorentz}.
\end{remark}

\subsection{An eigenvalue problem}
The P\'{o}lya-Szeg\"{o} inequality allows us to obtain a sharp lower bound 
for the first eigenvalue of the following nonlinear eigenvalue problem
\begin{equation} 
 \label{P2} 
\left\{\begin{array}{ll}
-{\rm div} (|\nabla u |^{p-2}  \nabla u ) =\lambda |x| ^{- \beta p} 
|u|^{p-2}u & \mbox{in }
\Omega
\\
u=0 & \mbox{on } \partial\Omega %\label{eq1}
\end{array}
\right.\\
\\
\end{equation}
where $  \Omega$ is a bounded domain in $\mathbb{R}^N$,  $1< p< N$
% $\beta p>N$ 
and  $0\leq \beta <1$.
This eigenvalue problem, together with some related elliptic problems for the $p$-Laplacian
has been studied in \cite{CP}. 

\noindent We set 

\begin{equation}
 \lambda_1 (\Omega ) = \min 
\left\{ 
\frac
{ \displaystyle\int_{\Omega}  | \nabla \varphi |^p  dx }
{ \displaystyle\int_{\Omega}  
|x|^{-\beta p }   
| \varphi (x) |^p  dx }: \
    \varphi \in  W^{1,p}_0 ( \Omega )\setminus \{ 0\}
	 \right\} \,.
\end{equation}
Observe that the bounds on $\beta$ and $p$ assure that the 
imbedding of $W^{1,p}_0(\Omega)$ in $L^{p}(\Omega, |x|^{-\beta p } dx)$ is compact 
(see, e.g. \cite{CP}).

The following result holds true

\begin{theorem}
\label{eigenvalue}
Let $\Omega ^{\star}$ denote the $\mu _{-\beta p} $-symmetrization of $\Omega $.
We have  
\begin{equation} \label{l1bis} 
\lambda_1(\Omega )\ge \lambda_1(\Omega ^{\star})
\end{equation}
where  $\lambda_1(\Omega ^{\star})$ is the first eigenvalue of the problem

\begin{equation}
\\
\\
\left\{\begin{array}{ll}
-{\rm div} (| \nabla v|^{p-2} | \nabla v| ) 
=
\lambda |x| ^{-\beta p}  |v|^{p-2}v & \mbox{ in }
\Omega^{\star}
\\
v=0 & \mbox{ on } \partial\Omega^{\star} .
\label{eq1}
\end{array}
\right.
\\
\end{equation}
\end{theorem}
   
\noindent 
{\sl Proof:\/} 
Put $l:= -\beta p$ and $k:= -\beta (p-1)$. Then it follows that $k\le 0$, $l+N >0$ and $l(N-1)/N -k \leq 0$. Hence the conditions 
%(\ref{k<l+1}) and 
(\ref{lk}) are satisfied. Furthermore, we have $pk+l(1-p)=0$. Applying Theorem \ref{ps}, we obtain, by the definition of $u^{\star}$,
\begin{eqnarray*}
 & & \int_{\Omega } |\nabla u|^p \, dx \geq \int_{\Omega ^{\star }} |\nabla u^{\star } |^p \,dx ,
\\
 & & 
\int_{\Omega }   |u(x)|^p |x|^{-\beta p } \, dx =\int_{\Omega^{\star }}    |u^{\star }(x)|^p|x|^{-\beta p}     \, dx, 
\end{eqnarray*}
and the result follows.
$\hfill \Box $   
\bigskip

\begin{remark} \rm 
Let $\Omega^\#$ be the ball centered at the origin  having the same Lebesgue measure 
as $\Omega$, that is, $\Omega ^\# $ is the Schwarz symmetrization ($=\mu _0$-symmetrization of $\Omega $). Then the following estimate holds
\begin{equation}
\label{l2}
 \lambda_1(\Omega )\ge\lambda_1(\Omega^\#) ,
\end{equation}
see  \cite{Bandle}. 
Indeed, if 
 $u^\# (x) $  denotes the Schwarz symmetrization of $u$, 
then the following estimate holds true
$$ 
\displaystyle 
\frac
{\displaystyle \int_\Omega  |\nabla u | ^p\, dx}
{\displaystyle \int_{\Omega }   |u(x)|^p |x|^{-\beta p} \, dx}\ge \displaystyle \frac{\displaystyle \int_{\Omega^\#} |\nabla u^\# |^p \, dx}{\displaystyle \int_{\Omega^\# }   | u^\#  (x)|^p |x|^{-\beta p }\, dx}\ge \lambda_1(\Omega^\#)
$$
which implies \eqref{l2}.

\noindent Observe that estimate (\ref{l2}) is worse than (\ref{l1bis}) since  $|\Omega^{\star}  |\le |\Omega^\#|$ implies
\begin{equation}
 \lambda_1(\Omega^\#)\le\lambda_1(\Omega^{\star}).
\end{equation}
\end{remark}


\begin{thebibliography}{99} 
\renewcommand{\baselinestretch}{0.9}
\small 

\bibitem{Alvino} 
{\sc  A. Alvino}, 
Sulla diseguaglianza di Sobolev in spazi di Lorentz, {\sl Bollettino U.M.I} {\bf 14-A} (1977), 148-156.

\bibitem{AH} 
{\sc  H. Ando, T. Horiuchi},
 On the weighted rearrangement of functions and degenerate nonlinear elliptic equations, {\sl Math. J. Ibaraki Univ.} {\bf 44} (2012), 17-31.

\bibitem{Bandle} 
{\sc  C. Bandle}, Isoperimetric inequalities and applications.
{\sl Monographs and Studies in Mathematics}, 7. Pitman (Advanced Publishing Program), Boston, Mass.-London, 1980. x+228 pp. 

\bibitem{BBMP} 
{\sc M.F. Betta, F. Brock, A. Mercaldo, M.R. Posteraro}, 
A weighted isoperimetric inequality and applications to symmetrization. 
{\sl J. Inequal. Appl.} {\bf 4} (1999), no. 3, 215–-240. 

\bibitem{BBMP2} 
{\sc M.F. Betta, F. Brock, A. Mercaldo, M.R. Posteraro}, Weighted isoperimetric inequalities on $\mathbb{R}^N $ and applications to rearrangements.
{\sl Math. Nachr.} {\bf 281} (2008), no. 4, 466--498.



\bibitem{BBCLT}   
{\sc W. Boyer, B. Brown, G. Chambers, A. Loving, S. Tammen}, 
Isoperimetric regions in $\mathbb{R}^n $ with density $r^p $,  
ArXiv:1504.01720v2.

\bibitem{BrasPhil} 
{\sc L. Brasco, G. De Philippis, B. Ruffini},  Spectral optimization for the Stekloff-Laplacian: the stability issue. {\sl J. Funct. Anal.} {\bf 262} (2012), no. 11, 4675–-4710.

\bibitem{BCM} 
{\sc F. Brock, F. Chiacchio, A. Mercaldo}, A class of
degenerate elliptic equations and a Dido's problem with respect to a
measure. {\sl J. Math. Anal. Appl.} {\bf 348} (2008), no. 1, 356--365.

%\bibitem{BCM2} 
%{\sc F. Brock, F. Chiacchio, A. Mercaldo}, Weighted
%isoperimetric inequalities in cones and applications. {\sl Nonlinear
%Analysis T.M.A.} {\bf 75} (2012), no. 15, 5737--5755.

%\bibitem{BCM3} 
%{\sc F. Brock, F. Chiacchio, A. Mercaldo}, 
%A weighted isoperimetric inequality in an orthant. 
%{\sl Potential Anal.} {\bf 41} (2012),   171--186.


\bibitem{BMP} 
{\sc F. Brock, A. Mercaldo, M.R. Posteraro}, On Schwarz and
Steiner symmetrization with respect to a measure. 
{\sl Revista Matem\'{a}tica Iberoamericana} {\bf 29} (2013), 665--690.

\bibitem{XR} 
{\sc X. Cabre, X. Ros-Oton}, Sobolev and isoperimetric
inequalities with monomial weights. 
%arXiv:1210.4487, 2012.
{\sl J. Differential Equations } {\bf 255}  (2013),  4312--4336. 

\bibitem{XRS} 
{\sc X. Cabre, X. Ros-Oton, J. Serra}, Euclidean balls
solve some isoperimetric problems with nonradial weights. 
%arXiv:1210.1788,2012.
{\sl  C. R. Math. Acad. Sci. Paris} {\bf 350}  (2012),  945--947

\bibitem{CKN} 
 {\sc L.  Caffarelli, R.  Kohn, L. Nirenberg},
 First order interpolation inequalities with weights. {\sl Compositio Math.} 
 {\bf 53} (1984), no. 3, 259–-275.

\bibitem{CaldMus}
{\sc P. Caldiroli, R. Musina}, Symmetry Breaking of Extremals for the
Caffarelli-Kohn-Nirenberg Inequalities
in a Non-Hilbertian Setting, {\sl Milan J. Math.} {\bf 81} (2013), 421–-430.
 
 \bibitem{CMV} 
{\sc A. Ca\~{n}ete, M. Miranda Jr., D. Vittone}, Some
isoperimetric problems in planes with density. {\sl J. Geom. Anal.} 
{\bf 20} (2010), no.2, 243--290.

 
 \bibitem{CJQW} 
{\sc T. Carroll, A. Jacob,  C. Quinn, R. Walters}, The isoperimetric problem on planes with density. {\sl Bull. Aust. Math. Soc.  }
{\bf 78}  (2008), no.2, 177--197.

 
 \bibitem{CRT} 
{\sc D. Cassani,  B. Ruf, C. Tarsi,}
Optimal Sobolev type inequalities in Lorentz spaces. {\sl Potential Anal. } {\bf 39} (2013), no. 3, 265--285. 

\bibitem{CW} 
{\sc F.  Catrina, Z. Wang},
On the Caffarelli-Kohn-Niremberg inequalities: sharp constants,
existence (and nonexistence),
{\sl Comm. Pure Appl. Math.} {\bf 54} (2001), no 2, 229--258.

\bibitem{Cham}
{\sc G. R. Chambers}, Proof of the Log-Convex Density Conjecture, 	arXiv:1311.4012v3 

\bibitem{ChiHo}
{\sc N. Chiba, T. Horiuchi}, On radial symmetry and its breaking in the Caffarelli-Kohn-Nirenberg inequalities for $p=1$. {\sl Math. J. Ibaraki Univ.} {\bf 47} (2015), 49--63. 

\bibitem{CP} 
{\sc E. Colorado, I. Peral}, Eigenvalues and bifurcation for elliptic equations with mixed Dirichlet-Neumann boundary conditions related to Caffarelli-Kohn-Nirenberg inequalities. 
{\sl Topological Methods in Nonlinear Analysis}, 
Journal of the Juliusz Schauder Center,
{\bf 23}, (2004), 239 --273.

\bibitem{C} 
{\sc  G. Csat\'o}, 
An isoperimetric problem with density and the Hardy Sobolev inequality in ${\mathbb R}^2$, {\sl Differential Integral Equations} {\bf 28} (2015), no. 9-10, 971–-988.

\bibitem{DDNT} 
{\sc  J. Dahlberg, A. Dubbs, E. Newkirk, H. Tran},
Isoperimetric regions in the plane with density $r^p$, 
{\sl New York J. Math.} {\bf 16} (2010), 31-51. 

\bibitem{DHHT} 
{\sc  A. Diaz, N. Harman, S. Howe, D. Thompson}
Isoperimetric problems in sectors with density. {\sl Adv. Geom. } {\bf 12} (2012), 589-619.




\bibitem{DolEstLoss}
{\sc J. Dolbeault, M. Esteban, M. Loss}, Rigidity versus symmetry breaking via nonlinear flows on cylinders and euclidean spaces, 
 (2015), ArXiv:1506.03664v1.
 
 
 \bibitem{E}  
{\sc  H. Eqnell}, 
Elliptic boundary value problems with singular coefficients and 
critical nonlinearities. {\sl Indiana Univ. Math. J.} {\bf 38} (1989), 235--251.

\bibitem{FleRi} 
{\sc W.H. Fleming, R. Rishel}, An integral formula for total gradient variation. 
{\sl Arch. Math. (Basel)} {\bf 11} (1960),  218–-222.  

\bibitem{G} 
{\sc E. Giusti }, Minimal surfaces and functions of bounded variation. {\sl Monographs
in Mathematics, 80. Birkh\"{a}user Verlag}, Basel, 1984.

\bibitem{Hardy} 
{\sc G.H. Hardy}, Notes on some points in the integral calculus. {\sl Messenger Math.} {\bf 48} (1919), 107–-112.

\bibitem{H} 
{\sc T. Horiuchi}, Best constant in weighted Sobolev inequality with weights being powers of distance from the origin. {\sl J. Inequal. Appl.} {\bf 1} (1997), no. 3, 275–-292.

\bibitem{HK} 
{\sc T. Horiuchi, P. Kumlin},
 On the Caffarelli-Kohn-Nirenberg-type inequalities 
involving critical and supercritical weights.
{\sl Kyoto J. Math.} {\bf 52} (2012), no. 4, 661–-742.

\bibitem{Howe}
{\sc S. Howe}, The Log-Convex Density Conjecture and vertical surface area in warped products. 
{\sl Adv. Geom.} {\bf 15} (2015), 455–-468.

\bibitem{Kaw}
{\sc B. Kawohl}, Rearrangements and convexity of level sets. Springer-Verlag N.Y. (1985). 

\bibitem{KZ} 
{\sc A.V. Kolesnikov, R.I. Zhdanov}, 
On isoperimetric sets of radially symmetric measures. Concentration, functional inequalities and isoperimetry, 123–-154,
{\sl Contemp. Math.} {\bf 545}, Amer. Math. Soc., Providence, RI, 2011.

\bibitem{Lan} 
{\sc R. Landes}, Some remarks on rearrangements and functionals with non-constant density. {\sl Math. Nachr.} {\bf 280}, no. 5-6, 560--570 (2007)

\bibitem{LL} 
{\sc N. Lam, G. Lu}, Sharp constants and optimizers for a class of the Caffarelli-Kohn-Niremberg inequalities, (2015) Arxiv:1510.01224v1.


\bibitem{M} 
{\sc V. Maz'ja}, Lectures on isoperimetric and isocapacitary inequalities in the theory of Sobolev spaces. Heat kernels and analysis on manifolds, graphs, and metric spaces (Paris, 2002), 307--340, 
{\sl Contemp. Math.} {\bf 338}, Amer. Math. Soc., Providence, RI, 2003. 

%\bibitem{Mo} {\sc F. Morgan}, The log-convex conjecture. 
%https://sites.williams.edu/Morgan/2010/04/03/the-log-convex-density-conjecture/. 
%Blog on the homepage of F. Morgan.

\bibitem{Mo} 
{\sc F. Morgan}, Manifolds with density. \textsl{Notices
Amer. Math. Soc.} \textbf{52} (2005), no.8, 853--858.

\bibitem{Mo2} 
{\sc F. Morgan}, The Log-Convex Density Conjecture. \textsl{%
Contemporary Mathematics} \textbf{545} (2011), 209--211.


\bibitem{Musina} 
{\sc R. Musina}, Weighted Sobolev spaces of radially symmetric functions. 
{\sl Ann. Mat. Pura Appl.} (4) {\bf 193} (2014), no. 6, 1629–-1659.

\bibitem{Stein} 
{\sc E. Stein}, 
Singular integrals and differentiability properties of functions. 
{\sl Princeton Mathematical Series}, no. 30, Princeton University Press, Princeton, N.J. 1970.

\bibitem{Talenti1} 
{\sc G. Talenti}, 
The standard isoperimetric theorem. {\sl Handbook of convex geometry}, Vol. A, B, 73–-123, North-Holland, Amsterdam, 1993. 

\end{thebibliography}
\end{document}